\documentstyle[twoside,12pt]{amsart}
\input psfig

\def\SBIMSMark#1#2#3{
 \font\SBF=cmss10 at 10 true pt
 \font\SBI=cmssi10 at 10 true pt
 \setbox0=\hbox{\SBF Stony Brook IMS Preprint \##1}
 \setbox2=\hbox to \wd0{\hfil \SBI #2}
 \setbox4=\hbox to \wd0{\hfil \SBI #3}
 \setbox6=\hbox to \wd0{\hss
             \vbox{\hsize=\wd0 \parskip=0pt \baselineskip=10 true pt
                   \copy0 \break%
                   \copy2 \break%
                   \copy4 \break}}
 \dimen0=\ht6   \advance\dimen0 by \vsize \advance\dimen0 by 8 true pt
 \dimen2=\hsize \advance\dimen2 by .25 true in
 \ht6=0pt \dp6=0pt
 \setbox8=\vbox to \dimen0{\vfill \hbox to \dimen2{\hss \copy6}}
 \ht8=0pt \dp8=0pt \wd8=0pt
 \copy8
}

\def\dist{\text{dist}}

\def\diam{\text{diam}}
\def\Haus{{\cal H}}
\theoremstyle{plain}                    
\newtheorem{thm}{Theorem}[section]
\newtheorem{cor}[thm]{Corollary}
\newtheorem{lemma}[thm]{Lemma}

\theoremstyle{definition}               

\newcounter{rem}[section]
\newenvironment{remark}{\addtocounter{rem}{1}
                        {\bf Remark \thesection .\therem:}}{}

\numberwithin{equation}{section}
\numberwithin{figure}{section}

\newcommand{\mycaption}[1]{\caption{\fontsize{10}{12}\selectfont #1}}
\begin{document}
\SBIMSMark{1998/6}{May 1998}{}

\renewcommand{\thepage}{}              
\title  [A non-removable set]
        {Non-removable sets  for quasiconformal  \\
         and locally biLipschitz  mappings in ${\Bbb R}^3$.}


\author {Christopher J. Bishop}
\address{C.J. Bishop\\
         Mathematics Department\\
         SUNY at Stony Brook \\
         Stony Brook, NY 11794-3651}
\email {bishop@@math.sunysb.edu}
\thanks{The  author is partially supported by NSF Grant DMS 95-00557 and
        an Alfred P. Sloan research fellowship.}


\begin{abstract}
We give an example of a totally disconnected set $E \subset 
{\Bbb R}^3$  which is not removable for quasiconformal homeomorphisms,
i.e., there is a homeomorphism $f$ of ${\Bbb R}^3$ to itself
which is quasiconformal off $E$, but not quasiconformal on 
all of ${\Bbb R}^3$. The set $E$ may be taken with Hausdorff 
dimension $2$.  The construction also gives a non-removable
set for locally biLipschitz  homeomorphisms.
\end{abstract}

\maketitle

%
\bigskip
\setcounter{page}{1}
\renewcommand{\thepage}{\arabic{page}}
\section{Statement of results}  \label{statement}

If a homeomorphism of ${\Bbb R}^d$ to itself is
 quasiconformal except on  a compact set
$E$, does it have to be quasiconformal on all of ${{\Bbb R}}^d$?
If so, $E$ is called removable for quasiconformal mappings.
The purpose of this paper is to construct examples of non-removable
sets in ${\Bbb R}^3$
 which are  as small as possible, both  topologically  (they are 
totally disconnected) and   metrically (they have Hausdorff dimension 
$2$). 

A mapping is called {\it quasiconformal } on $\Omega \subset 
{\Bbb R}^d$ if there is an $M< \infty$
 so that 
$$ \limsup_{r \to 0} \frac { \sup_{|x-y|=r} |f(x)-f(y)|}{\inf_{|x-y|=r}
                  |f(x)-f(y)|} \leq M \quad \forall x \in \Omega .$$
(See \cite{Heinonen-Koskela95} or Theorem 34.1 of  \cite{Vaisala71}.)
Our method will actually give non-removable sets for an even more 
restrictive class of mappings. We say that a mapping is 
{ \it locally biLipschitz}  on $\Omega$ if there is an 
$M < \infty$ so that  for every $x \in \Omega$ there is an $r=r(x)>0$ so that
$|x-y| < r$ implies 
$$ M^{-1} \leq \frac {|f(x)-f(y)|}{|x-y|} \leq M.$$
Such mappings are also called {\it  bounded length distortion} 
(e.g., \cite{Vaisala87}, \cite{Vaisala89}) or {\it local quasi-isometries}
(e.g., \cite{Gehring82}, \cite{John68}).
If a quasiconformal mapping is biLipschitz on  dense open set 
then it is globally biLipschitz, and hence 
a non-removable set for the biLipschtiz maps 
  is also non-removable for quasiconformal maps.

\begin{thm} \label{main}
There is a totally disconnected set $E \subset {\Bbb R}^3$ which 
is nonremovable for 
locally biLipschitz (and hence for quasiconformal) maps.
 If $\varphi(t) = o(t^2)$ then we may choose $E$ and $f$ so that 
$\Haus^\varphi(E) = \Haus^\varphi(f(E))=0$.
\end{thm}

Here $\Haus^\varphi$ denotes the $\varphi$-Hausdorff measure, i.e.,
$$ \Haus^\varphi(E) = \lim_{\delta\to 0} [ \inf\{ \sum_j \varphi(r_j), 
E \subset \cup_j B(x_j, r_j), r_j \leq \delta \}].$$
Our result is sharp in the sense that if $\Haus^\varphi(E) =0$ for every
$\varphi$ such that $\varphi(t) = o(t^2)$, then $E$ has $\sigma$-finite
$\Haus^2$ measure (\cite{Besicovitch56}) and  hence is removable for 
homeomorphisms which are quasiconformal off $E$
 (Theorem 35.1 of \cite{Vaisala71}).
Since locally biLipschitz mappings have gradient in $L^\infty$ 
on $\Omega$ we see that our examples are also non-removable 
for the Sobolev spaces $W^1_p$ for every $p \leq \infty$, 
answering a question of P. Koskela.

The only previously known examples of nonremovable sets in ${\Bbb R}^3$
  either have interior (trivial) or 
are of the form $F= E\times [0,1]^2$ for any uncountable $E \subset {\Bbb R}$ 
(\cite{Carleson50}, \cite{Kaufman-Wu96}). 
In the latter case, assume  $E\subset [0,1]$. It  supports 
a non-atomic probability measure $\mu$ which is singular to Lebesgue 
measure. If we define $f$ to be the identity 
outside $ S= \{ (x,y) \in {\Bbb R} \times {\Bbb R}^2: 0<x<1, |y| < 1\}$
 and 
$$ f(x,y) = (|y|x+(1-|y|)(\frac {x+\mu([0,x])}{2}), y),$$
inside $S$
then we easily see that $f$ is  a homeomorphism of ${\Bbb R}^3$
which is  locally biLipschitz on the 
complement of $F$, but maps a set of zero volume  to positive volume,
and hence is not even quasiconformal on all of ${\Bbb R}^3$.
 See \cite{Bishop94b}, \cite{Gehring60}, \cite{Kaufman84} and  \cite{Wu97}
for other  constructions of non-removable sets in ${\Bbb R}^2$.


One of the most striking aspects of the construction is that it allows
 one to approximate
any smooth diffeomorphism by quasiconformal or locally biLipschitz
 maps with uniform bounds on the constants (independent of the map
 being approximated), as long as we  ``throw out''  a fairly small set.  

\begin{cor} \label{app-cor}
Suppose $\Omega_1$ and $\Omega_2$  are open sets in ${\Bbb R}^n$ which are
diffeomorphic.
Then for any $\epsilon >0$, 
 there is a homeomorphism $f: \Omega_1 \to \Omega_2$ which
is  quasiconformal  except on a totally disconnected set
 $E$ and which approximates $h$ to within $\epsilon$. 
 For any measure function $\varphi(t) = o(t^2)$ we may take
 $\Haus^\varphi (E) = \Haus^\varphi(f(E))=0$. If 
$\Omega_1$ and $\Omega_2$ are diffeomorphic by a volume preserving 
map we may take $f$ to be locally biLipschitz except on $E$.
\end{cor}

For conformal mappings in the plane, this type of result was proved 
by Gehring and Martio (Theorem 4.1, \cite{Gehring-Martio85}). They 
showed  that
there exists a Cantor set $E$ in the unit disk ${\Bbb D}$ so that ${\Bbb D} 
\setminus E$ is conformally equivalent to the plane minus a Cantor set.
The two dimensional version of our construction gives a geometric 
construction of  such a set. In ${\Bbb R}^3$ it shows that there is a
Cantor set $E$ in the unit ball ${\Bbb B}$ so that ${\Bbb B} \setminus E$
is quasiconformally equivalent to ${\Bbb R}^3$ minus a Cantor set.
(${\Bbb R}^3$ is not itself quasiconformally equivalent to ${\Bbb B}^3$, 
e.g., Section 17.4 of  \cite{Vaisala71}.)

By a result of Dacorogna and Moser \cite{Dacorogna-Moser}, if 
$\Omega_1$ and $\Omega_2$  have smooth boundaries, 
are diffeomorphic up to the boundary 
and have the same  volume, then there is a diffeomorphism 
between them which preserves volumes. Thus the last statement 
in Corollary \ref{app-cor} is fairly general.

Suppose  $E$ is our non-removable set of dimension $2$ and 
$f$ is quasiconformal on all of ${\Bbb R}^3$. Then $f(E)$ 
must also be non-removable hence have dimension $\geq 2$.
On the other hand, our construction of $E$ will show that it
is ``tame'', i.e.,  it can 
be mapped to the standard Cantor set by a homeomorphism of ${\Bbb R}^3$
(and hence it can be mapped to a set of dimension zero by some
homeomorphism).  
Let $H({\Bbb R}^d)$ be the collection of homeomorphisms of ${\Bbb R}^d$
to itself and $QCH({\Bbb R}^d)$ be the subset of quasiconformal 
homeomorphisms.  For $E \subset {\Bbb R}^d$, define
$$\dim_{H}(E) = \inf_{f \in H({\Bbb R}^d)} \dim(f(E)),$$
$$\dim_{QC}(E) = \inf_{f \in QCH({\Bbb R}^d)} \dim(f(E)).$$

\begin{cor}
There is a compact $E \subset {\Bbb R}^3$  
such that $\dim(E) = \dim_{QC}(E) =2$, but $\dim_H(E)=0$.
\end{cor}

There are at least two other types of Cantor set whose dimension can't 
be lowered by quasiconformal maps. First, 
since quasiconformal maps are absolutely continuous
with respect to Lebesgue measure, a Cantor set of positive measure
has this property. More generally, if $ E \subset {\Bbb R}^d$
and  $\dim(E) = d$ then $\dim_{QC}(E) = d$ by results of 
\cite{Gehring-Vaisala71}.
Second, there are totally disconnected  sets $F$ (e.g., Antoine's necklace,
\cite{Antoine}, \cite{Hocking-Young}) 
whose complement is not simply connected, and hence 
$\dim_{H}(E) \geq 1$.
For any $0<\alpha<d$ is there a compact  $E \subset {{\Bbb R}}^d$
with  $ \dim_{QC}(E) = \alpha  = \dim(E)$? The only known
examples are when $\alpha$ is an integer. Is $\dim_{H}(E)$ always
an integer?

An open set,  $\Omega = {\Bbb R}^n \setminus E$, is called 
{\it quasiconvex} if there is a $C < \infty$ such that  any two points in 
$\Omega$ can be joined by a path in $\Omega$ of length at most 
$C|x-y|$.  If $\Omega$ is quasiconvex and $E$ has zero 
measure,  then $E$ must be 
removable for locally biLipschitz maps  (Lemma \ref{Q-convex}).
Our construction can be modified to
give a non-removable set $E$  for quasiconformal mappings whose complement
is quasiconvex and hence is removable for locally biLipschitz 
maps. Thus the two classes of sets are distinct. 

In fact, we can 
considerable strengthen the quasiconvexity as follows.
In the terminology of \cite{Kaufman-Wu96}, $E$ is called a {\it
weak porus set} if each $x \in E$ is contained in  a sequence of 
cubes $Q_j$ with diameters tending to zero and such that
 $(Q_j \setminus (1-\alpha_j)Q_j ) \cap E  = \emptyset$ for some 
positive sequence $\{ \alpha_j\}$. This property implies 
$E$ is a totally disconnected set  in a strong way,
and easily implies the complement is quasiconvex.



\begin{cor} \label{porus-cor} 
There is a weak porus set $ E \subset {\Bbb R}^3$ which is non-removable
for quasiconformal mappings. 
\end{cor}

As in Theorem \ref{main} we may take $\Haus^\varphi(E) =0$. 
Although $\partial Q_j$ misses $E$, it must be very  close to $E$
in the following sense. A result of Kaufman and Wu \cite{Kaufman-Wu96}
says that if $E$ is weakly porus with sequence $\{\alpha_j\}$ and
$\sum \alpha_j = \infty$, then $E$ is removable for quasiconformal
mappings (this generalizes  a result of 
 Heinonen and Koskela \cite{Heinonen-Koskela95} with 
$\alpha_j = \alpha$ independent of $j$).
Thus in our example, $\sum \alpha_j < \infty$ (in fact, 
$\alpha_j \to 0 $ very fast).

Our construction of the weakly porus non-removable set $E$ actually shows 
that it is a subset of a product set, i.e., 

\begin{cor} \label{product-cor}
There is a Cantor set $E \subset [0,1]$ so that $E \times E \times E$
is non-removable for quasiconformal mappings in ${\Bbb R}^3$.
\end{cor}

Ahlfors and Beurling proved that a product set in the plane  
is removable if both  factors have zero length
(Theorem 10, \cite{Ahlfors-Beurling50}).  Is this 
is true for triple products in ${\Bbb R}^3$?
As noted earlier products of the form $E \times [0,1]^2$, $E \subset
{\Bbb R}$ are removable iff $E$ is countable. When are products 
$E \times [0,1]$, $E \subset {\Bbb R}^2$ removable?
Every set of positive area in ${\Bbb R}^2$ is non-removable for 
quasiconformal mappings (e.g., \cite{David87} or \cite{Kaufman-Wu96}).
Is  every set of positive volume in ${\Bbb R}^3$ non-removable?

I thank Juha Heinonen and Jang-Mei Wu for telling me about 
the problem of  constructing a totally disconnected non-removable set for
quasiconformal mappings.
I also thank Pekka Koskela, Aimo Hinkkanen and Seppo Rickman for 
listening to an early version of the construction and 
encouraging me to write it down. Similarly for the participants
in the March 1996 Oberwolfach meeting on function theory, whose
comments improved the exposition and suggested some of the 
corollaries discussed above.  I am grateful to Richard Stong
for explaining the annulus conjecture to me and its connection to 
the construction.

The rest of the paper is organized as follows.
\begin{description}
\item [Section 2] We build a non-removable set for locally 
biLipschitz maps in ${\Bbb R}^2$.
\item [Section 3]
We build a ``flexible square'' which is the main building block 
of the three dimensional construction.
\item [Section 4]  We give the 
construction for locally biLipschitz mappings in ${\Bbb R}^3$.
\item [Section 5]   We show how to build non-removable sets  quasiconformal
maps  in the plane so that both  $E$ and $f(E)$ are  small. 
\item [Section 6] We modify the previous section to work in ${\Bbb R}^3$.
\item [Section 7] We prove Corollaries \ref{porus-cor} and \ref{product-cor}.
\item [Section 8] We construct non-removable sets for locally 
        biLipschitz maps with the additional property 
        that $\Haus^\varphi(E)=0$.
\item [Section 9]  We show how to get $\Haus^\varphi(f(E)) =0$.
\end{description}

\section{A non-removable set for locally biLipschitz mappings in ${\Bbb R}^2$}
\label{biLip-plane}

It is clear that an arbitrary smooth  mapping  $[0,1] \to {\Bbb R}^2 $ cannot
be approximated by a biLipschitz mapping with a uniform 
constant.
However, it can be approximated by a locally biLipschitz map 
if the line segment is replaced by an appropriately 
``wild'' arc. More precisely, 

\begin{lemma} \label{bad-arc-lem}
Suppose $g$ is a smooth homeomorphism from a neighborhood  $U$ of 
$[0,1]$  to ${\Bbb R}^2$.
For any $\epsilon>0$ there is an arc $\gamma$ (depending 
on $g$ and $\epsilon$) 
with the following properties.
\begin{enumerate}
\item $\gamma$ has endpoints $0$ and $1$.
\item $\gamma \subset [0,1]\times [-\epsilon , \epsilon] \subset U$.
(i.e., it approximates $[0,1]$ in the Hausdorff metric.)
\item 
      There is a locally biLipschitz map $f$  defined
      on a neighborhood of $\gamma$, so that
      for all $z \in \gamma$, $ | f(z) - g(z)| \leq \epsilon.$
\end{enumerate}
\end{lemma}

\begin{pf}
Consider the arc illustrated 
in Figure \ref{stretch-arc}. Although it is drawn a polygonal 
arc for simplicity, one should think of it as smooth (just round
the corners). Depending 
on the height, width and number of oscillations the 
arc can be stretched as much as we wish, by a length preserving 
map of the arc. The map can be 
extended to be locally biLipschitz in a neighborhood of the arc.
By taking an intermediate 
version of the arc, we obtain an
 arc which can be either stretched or contracted.
\begin{figure}[htbp]
 \centerline{ \psfig{figure=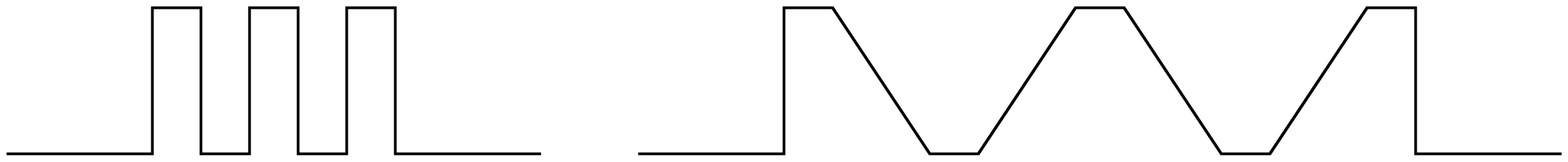,height=.5in} }
 \mycaption{ \label{stretch-arc} A stretchable arc } 
\end{figure}
\begin{figure}[htbp]
\centerline{ \psfig{figure=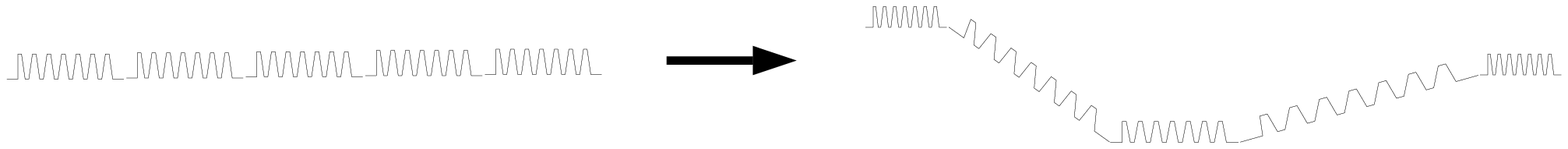,height=.5in} }
 \mycaption{ \label{spring-arc} Approximating polygonal arc with biLipschitz
images } 
\end{figure}
Using this building block and approximating by polygonal arcs
it is easy to see that we can approximate any smooth function 
on an appropriate $\gamma$.
\end{pf}

We will define our exceptional set $E$ as a limit
of sets $\{E^n\}$, each of which is a finite union of smooth
curves $E^n = \cup_j E^n_j$
 (with diameters tending to $0$ as $ n \to \infty$).
 Throughout the paper we
will label sets in the construction in the form $A^n_j$ where
the superscript $n$ denotes the generation of the construction and 
the subscript $j$ is an index enumerating  components in that
generation.

To begin the induction we start with $E^0 = \{z: |z|=1\}$.
 Let $\Omega^0$ be the complement of $E^0$,   let $\Omega^0_0$ 
be its unbounded component and $ \Omega^0_1$ the 
bounded component.
Define 
$$ f^0(z) = z, \quad z \in \Omega^0_0,$$
and extend $f^0$ to a smooth diffeomorphism (which we also call 
$f^0$) of the plane in anyway
you want (say with very large derivative at the origin).

The induction hypothesis is as follows.
Suppose we are given a compact set $E^n$
 which is a finite
union of $ J_n$ smooth closed curves, $\{ E^n_j\}$, which are 
disjoint with disjoint interiors. 
Let $\Omega^n$ be the complement
of $E^n$. Its  unbounded component is denoted $\Omega^n_0$  and 
the bounded components are denoted $\Omega^n_j$, $j=1,\dots,J_n$.
Suppose we are given a diffeomorphism  $f^n$ of the plane which 
is locally biLipschitz  on $\Omega^n_0$.
Let $Y^n_j$ be the bounded complementary component of $f^n(E^n_j)$.
Assume the  diameters  of $E^n_j$ and  $Y^n_j$ are less than $2^{-n}$.

\begin{figure}[htbp]
\centerline{ \psfig{figure=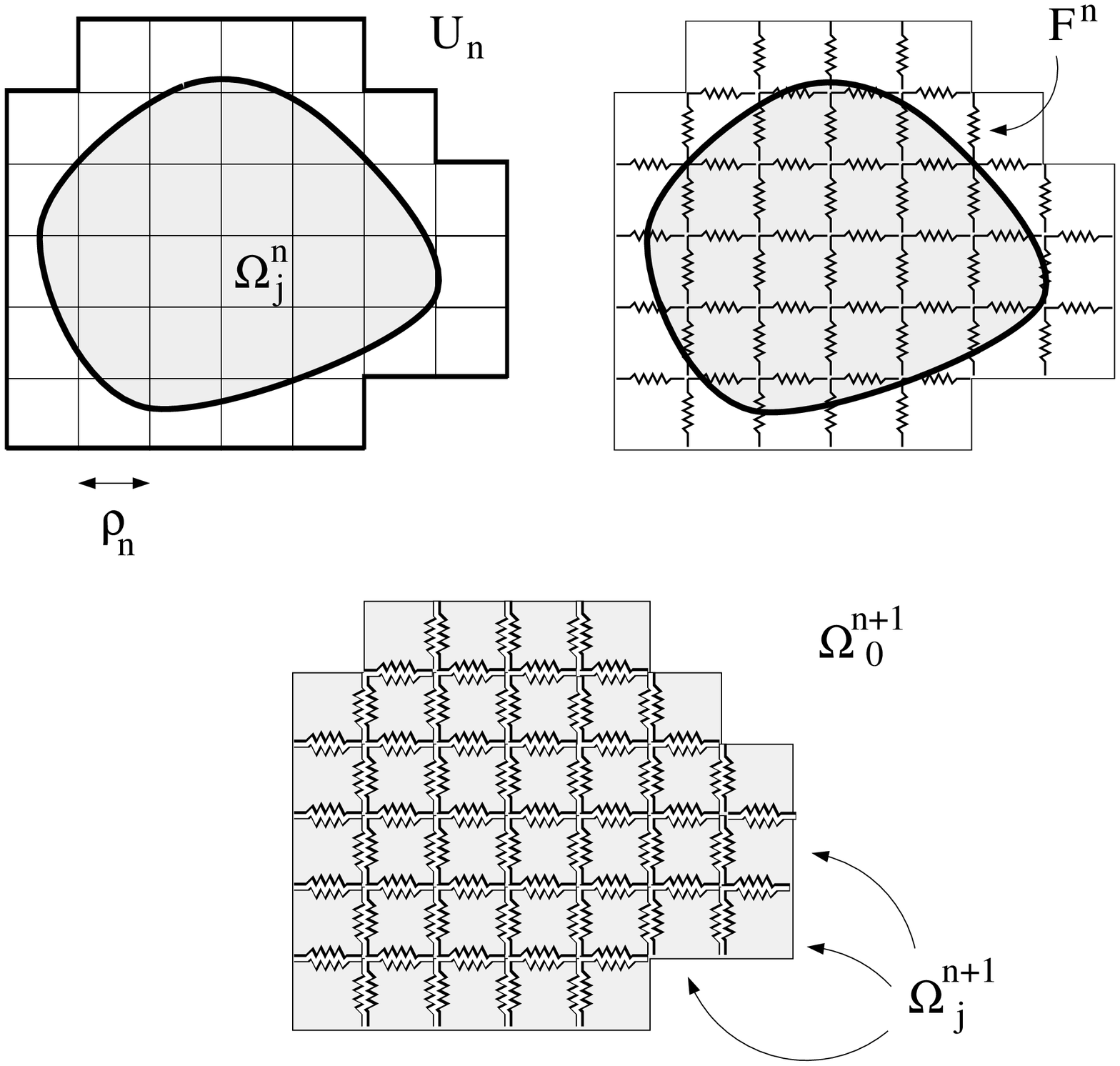,height=4in} }
 \mycaption{ \label{rough-cover} Covering $\Omega^n_j$ with squares}
\end{figure}

We now construct $f^{n+1}$
and $E^{n+1}$ from $f^n$ and $E^{n}$.
Choose a number $\rho_n \leq \frac 1{100} \min_{j \ne k} \dist(
E^n_j , E^n_k)$, and let ${\cal S}_n$ be a covering of 
of a neighborhood of $\overline{ \cup \Omega^n_j}$ by squares from the grid 
$\rho_n{{\Bbb Z}} \times  \rho_n{{\Bbb Z}}$. Let 
$U^n = \cup_{{\cal S}_n} Q$.
 By making $\rho_n$ even smaller,
if necessary, we may assume that $\diam(f^n(S)) \leq 2^{-n}/10$ for 
any of the squares $S$ in our cover.
 See Figure \ref{rough-cover} (upper left).

Replace each  interior 
edge of the union of squares
 by an arc from Lemma \ref{bad-arc-lem} and call the union 
of the arcs $F^n$. See Figure \ref{rough-cover} (upper right).
The arcs are chosen so that 
we can approximate $f^n$ by a locally biLipschitz map on 
a neighborhood of $F^n$. Let $\Omega^{n+1}_0$ be the union 
of $\Omega^n_0$ and this neighborhood. Without loss of generality, we
may assume $\Omega^{n+1}_0$ is bounded by a  finite number of smooth 
closed curves $\{E^{n+1}_j\}$ and that we have a $2$-biLipschitz 
homeomorphism $f^{n+1}$  defined on $\Omega^{n+1}$ which agrees 
with $f^{n}$ on $\Omega^n_0 \setminus U^n $.
See Figure \ref{rough-cover} (bottom).
 Define $f^{n+1}$ to be a diffeomorphism of
the plane by extending $f^{n+1}$ to the bounded complementary components
of $E^{n+1} = \partial \Omega^{n+1}_0$ in any way you want.
Finally, note that 
$$\diam(f^{n+1}(E^{n+1}_j)) \leq 
2 \cdot \diam( f^{n}(S)) \leq  2^{-n-1},$$
which is the final part of the induction hypothesis.

This completes the inductive step of the construction, i.e. given the
 set $E^n$ and mapping $\{f^n\}$ we have constructed $E^{n+1}$ and 
$\{ f^{n+1}\}$ which satisfy the induction hypothesis.
We now apply the following elementary lemma to the sets
$F^n = \overline{U^n}$.

\begin{lemma} \label{nested-sets}
Suppose $\{F^n \} \subset {\Bbb R}^d$  is decreasing, nested sequence of compact 
sets with disjoint  components $F^n = \cup_j F^n_j $, and 
$ \lim_{n\to 0} \sup_j \text{diam}(F^n_j) =0.$
Then $F= \cap_n F_n$ is totally disconnected. Suppose $\{ g^n\}$
is a sequence of homeomorphisms of ${\Bbb R}^d$ to itself such that
$ g^n = g^{n+1}$ on ${\Bbb R}^d \setminus F^n$ and 
$ \lim_{n\to 0} \sup_j \text{diam}((g^n(F^n_j)) =0.$
Then $\{ g_n\}$ converges uniformly to a homeomorphism 
$g: {\Bbb R}^d \to {\Bbb R}^d$.
\end{lemma}

We leave the proof of this to the reader.
Using the lemma we see our maps converge to a homeomorphism
which is locally biLipschitz off a totally disconnected set $E$.
Finally, to see that $g$ is not locally biLipschitz on all of ${\Bbb R}^2$, 
there are several things we could do. The easiest is to define the
homeomorphisms $\{f^n\}$  at each stage so that the limiting 
homeomorphism  is not H{\"o}lder of any positive order.
Thus it is not even quasiconformal on ${\Bbb R}^2$.


\section{A flexible square}
\label{flex-square-sec}

To do the construction in three dimensions, we follow the 
previous construction. However, when we get to the step where we replaced
each edge of the covering squares by a flexible arc, 
we will have to replace faces of a cube by flexible surfaces. 
Building such surfaces is the only difficult point in extending
the construction to higher dimensions.


\begin{lemma} \label{flex-square}
Suppose $g$ is a diffeomorphism 
of  a neighborhood of  $[0,1]^2$ into  ${\Bbb R}^3$. 
Suppose that $\epsilon >0$ is given. Then 
there is a surface $S_1$ and smooth mapping $G$ defined on 
a neighborhood of $S_1$ so that 
\begin{enumerate}
\item $S_1$ is a topological disk which approximates $[0,1]^2$ to
     within $\epsilon$ in the Hausdorff metric.
\item If $g$ is locally $M$-biLipschitz 
         on a neighborhood of $\partial [0,1]^2$ then 
       $G$ is locally  $(M + \epsilon)$-biLipschitz on $S_1$ and is 
       uniformly  locally biLipschitz   outside an 
        $\epsilon$-neighborhood of $\partial S_1$.
\item If $g$ is $M$-quasiconformal 
         on a neighborhood of $\partial [0,1]^2$ then 
       $G$ is $(M+\epsilon)$-quasiconformal on $S_1$.
\end{enumerate}
\end{lemma}

\begin{pf}
The basic idea is that the flexible surface can be obtained 
by ``folding'' a large square to make it oscillate, first in 
one direction and then in the other. We first show how to 
build a surface on which linear maps can be approximated.

Let $\gamma$ be the flexible arc constructed earlier and let 
$S_0$ be the surface obtained by crossing it with an interval.
See Figure \ref{easy-sheet2}.
\begin{figure}[htbp]
\centerline{ \psfig{figure=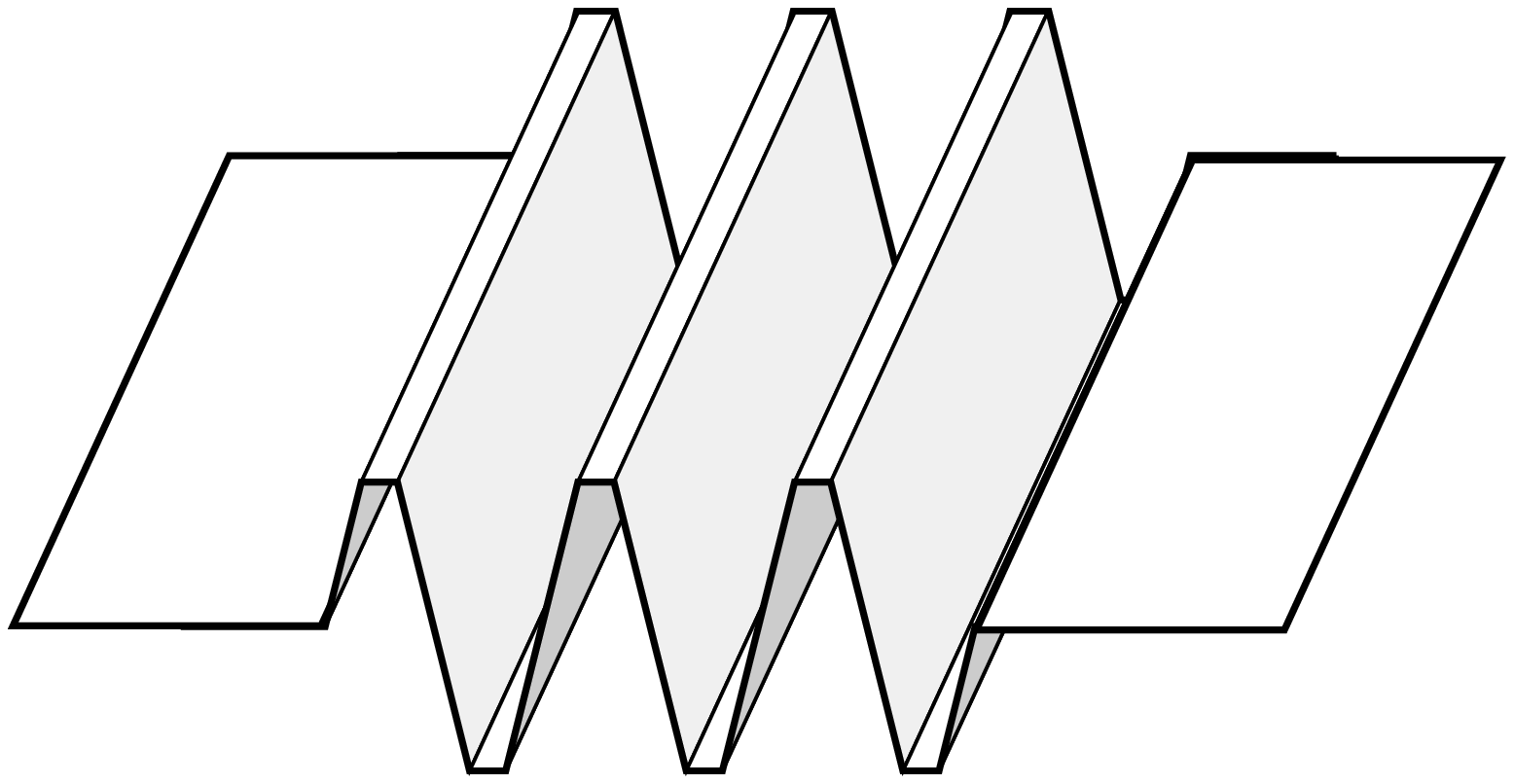,height=1.25in} }
 \mycaption{ \label{easy-sheet2} Surface flexible in one direction } 
\end{figure}
Then $S_0$ is can be stretched (by a locally biLipschitz map) in the 
direction parallel to $\gamma$ and may be skewed in the 
the perpendicular direction. Although the figure seems to 
have sharp corners, one should think of this as a smooth 
surface on small scales (or as polygonal with very small angle 
between adjacent faces).

Next, tile $S_0$ with small squares 
and replace each by a copy of $S_0$, but now with the 
copy of $\gamma$ oriented in the perpendicular direction.
The scale is chosen to be so small that $S_0$ looks flat 
in the small squares 
and so that adjacent tiles meet at very small angle. Thus adjacent 
tiles can be joined with only a small distortion.
Since $S_0$ can be stretched or shrunk in one direction and 
the small tiles can be stretched or shrunk in the other, the 
resulting surface $S_1$, can be simultaneously stretched 
or shrunk in both directions by a locally biLipschitz
map, i.e., we can approximate maps of the form 
$(x,y) \to (ax, by)$.
Drawing the surface itself is a bit complicated, but
Figure \ref{flex-surf-bend} gives an idea of what it looks like.
The picture is a little misleading because  the oscillations in 
different directions should be at very different scales.

\begin{figure}[htbp]
\centerline{ \psfig{figure=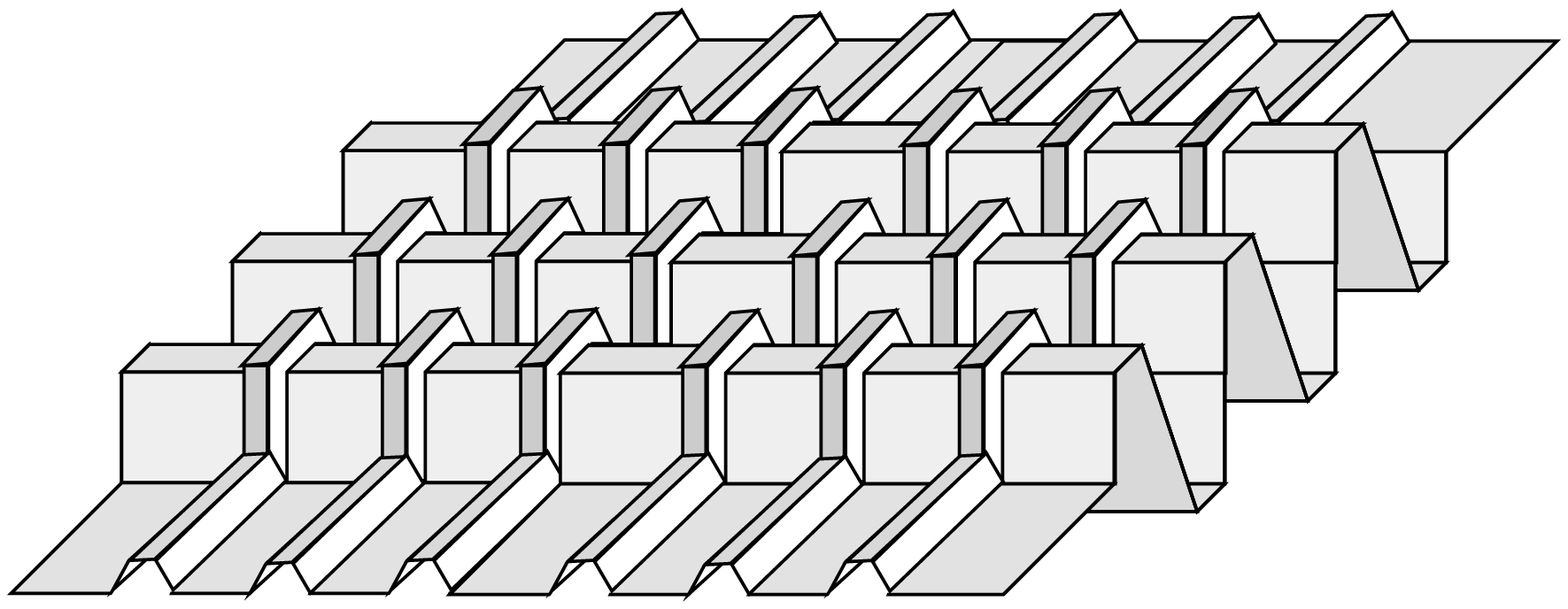,height=1.0in} }
\mycaption{ \label{flex-surf-bend} A flexible surface } 
\end{figure}

We may also assume that if the ``height'' of the large 
oscillations is $\delta$ then there is a $\delta \times
\delta$ flat square in each corner of $S_1$ and a $\delta$
wide strip along each edge in which there are only oscillations
in one direction. These strips will be used below to interpolate 
maps defined on adjacent squares.

The horizontal or vertical stretching is easy to see. The fact that 
we can skew our surface (i.e., approximate
maps of the form $(x,y) \to (x,y+cx)$) by a locally biLipschitz
map is a little harder but is illustrated in Figure \ref{skew-sqr}.
The picture shows the surface in Figure \ref{easy-sheet2} viewed 
from above. The white rectangles correspond to the horizontal 
pieces and the shaded rectangles to the almost vertical pieces.
First we stretch the square by making the almost vertical sides
horizontal. Then we apply a bounded distortion skew to each shaded
rectangle. Finally, we shrink in the horizontal direction by making
the shaded pieces almost vertical again. 

\begin{figure}[htbp]
\centerline{ \psfig{figure=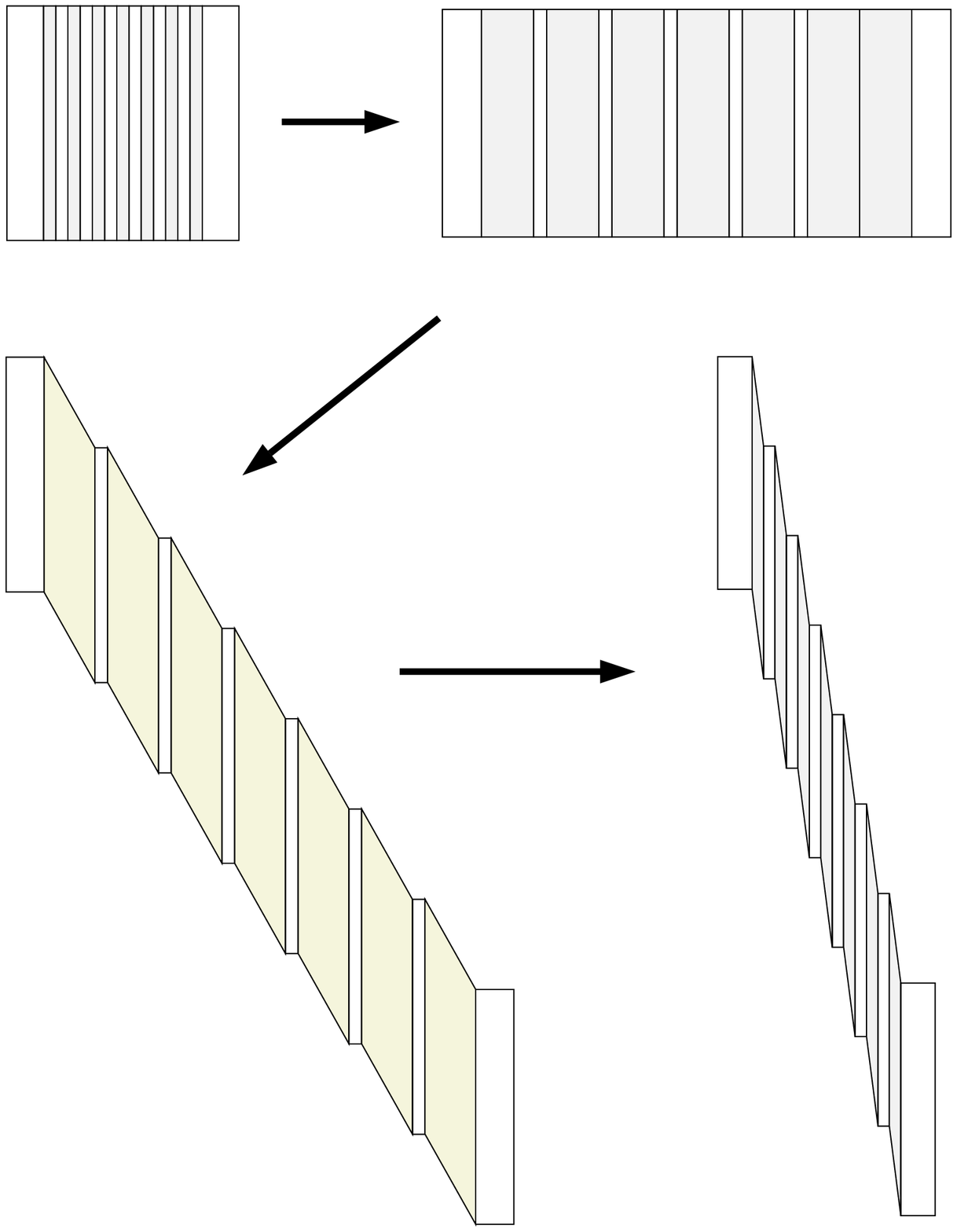,height=3in} }
\mycaption{ \label{skew-sqr} Approximating a skew map on the flexible surface} 
\end{figure}


Given any linear map of a square into the plane, 
we now have a surface which approximates the square and
a locally biLipschitz map which  approximates the  linear map.
Moreover, the degree of approximation is controlled in terms of 
the size of the oscillations of the surface. In particular, 
the images of the boundary arcs are, up to small distortion, 
simply the arcs stretched (or shrunk) to the appropriate diameter, 
i.e., up to a small distortion, the shape of the image arc is 
determined by the distance between its endpoints.
This is the main point which is used to glue together our approximations 
on adjacent flexible squares.

We now replace our smooth map $g$ by a piecewise linear 
approximation. Consider Figure \ref{lin-app}. It shows 
three regions; $U_1$, a neighborhood of $\partial [0,1]^2$, 
(the light gray region),  a square $Q =[\eta, 1-\eta]^2$ (the 
dark gray area), and  $U_2$, an open set which connects
the two. Inside $U_1$, we leave $g$ alone. We triangulate 
$U_2$ and define an approximation to $g$ which agrees with $g$ 
at the vertices on the of triangulation and on the faces 
which meet $U_1$, but which is linear on the faces which 
hit  $Q$. Since $g$ is smooth, we can do this and get a
locally biLipschitz approximation with constant close to that
for $g$ if we take the neighborhoods small enough.
\begin{figure}[htbp]
\centerline{ \psfig{figure=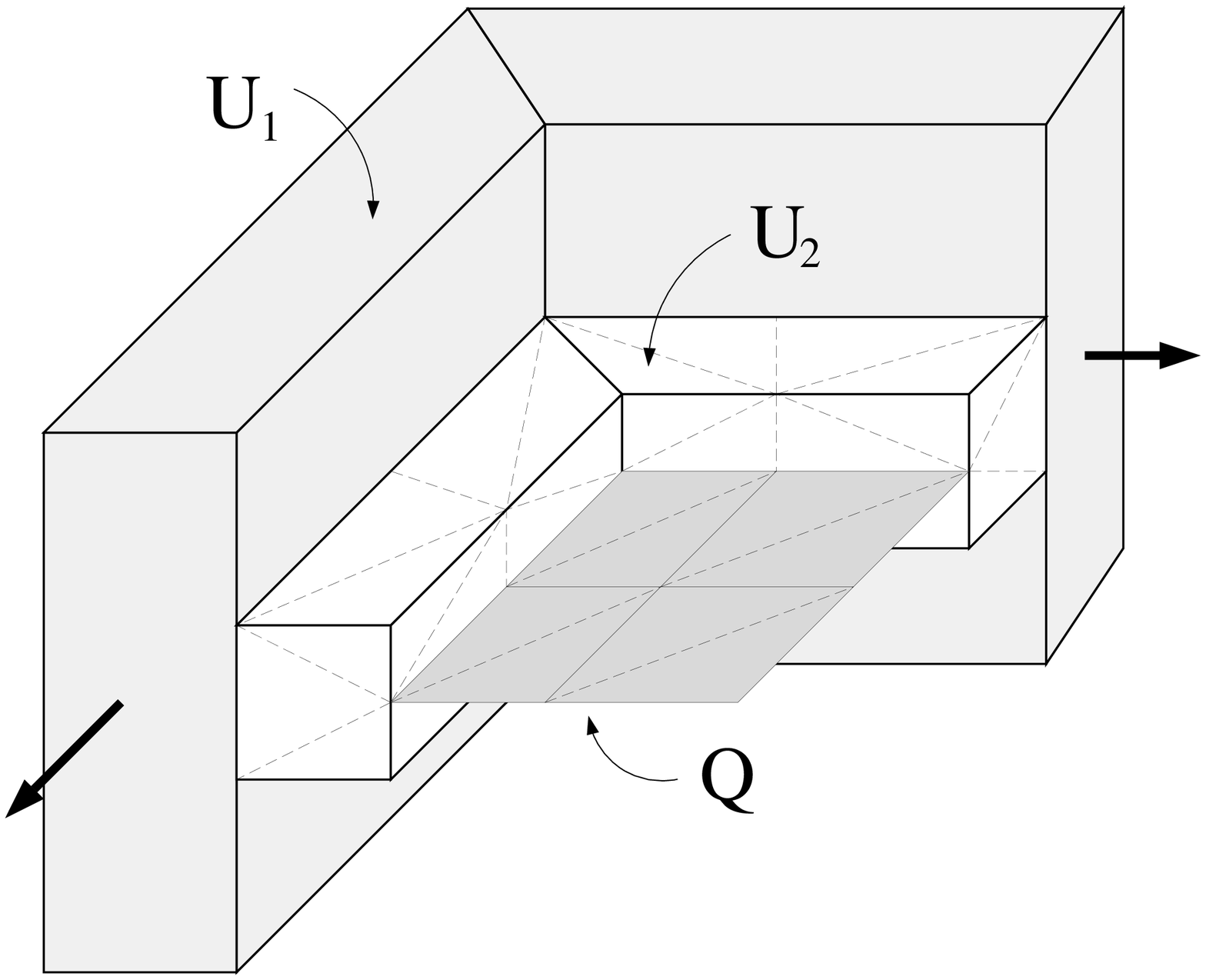,height=2.0in} }
 \mycaption{ \label{lin-app} Replace $g$ by a piecewise linear approximation}
\end{figure}

Divide $Q$ into small squares and 
divide each into two triangles by cutting it by 
a diagonal.
If the squares are small enough, then 
we  can replace $g$ by an approximation that agrees with $g$ 
at all the vertices  and which is   linear on each of
the triangles.

Now replace each subsquare in $Q$  with a copy of the same flexible
surface, chosen so that on each square we can approximate $g$ 
by a locally biLipschitz map on the surface and so that our 
approximation agrees with $g$ at the vertices of the triangulation
(First approximates the linear map on one of the two triangles 
whose union is the square; then apply a biLipschitz map with 
small distortion to ``bend'' the image along the diagonal to get
the fourth corner to agree.)

The maps defined on adjacent squares might not match up along the common 
boundaries but we can fix this as follows.
Along each edge of our flexible squares, we have a strip whose
width is greater than the 
vertical size of the oscillations and in which there are 
only oscillations in one direction. See Figure \ref{box-st-box}.
On the surface minus these
boundary strips we simply take the restriction of the map defined above. 
Inside the strip, the surface is a 
union of rectangles and in the image, the opposite sides are 
perturbed by a small angle and translation (this is due to 
our earlier remarks that the images of the boundary arcs are 
determined up to small distortion by the positions of the 
endpoints. The two boundary arcs we are trying to glue 
have the same endpoints and hence are small distortions of 
each other.) We can divide each rectangle into two triangles and
linearly interpolate the maps on the boundary arcs.
See Figure \ref{glue-edge}.
Similarly for the flat squares in the corners where four flexible 
squares come together.
The interpolated mapping is locally biLipschitz with a 
uniform bound.
\begin{figure}[htbp]
\centerline{ \psfig{figure=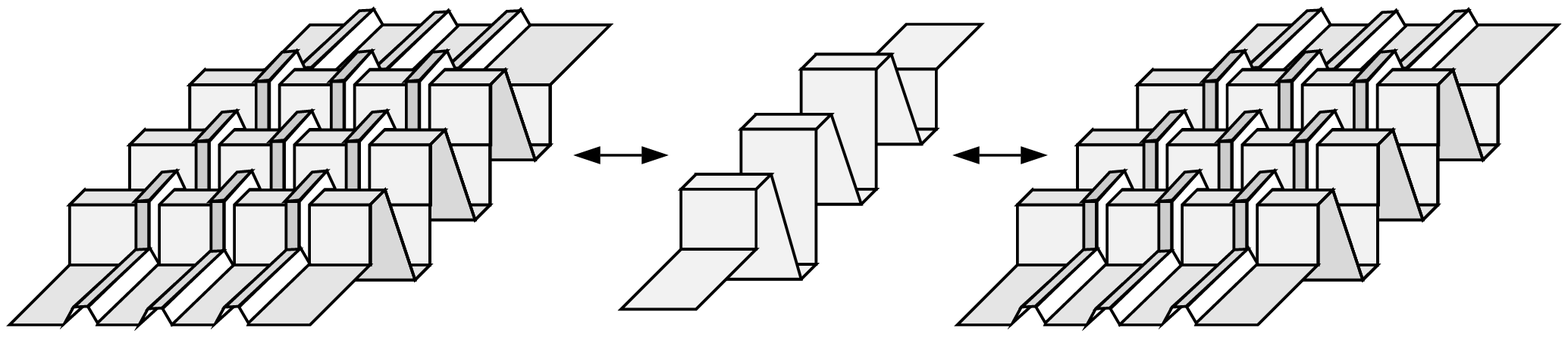,height=1.0in} }
 \mycaption{ \label{box-st-box} Strip adjoining adjacent squares } 
\end{figure}
\begin{figure}[htbp]
\centerline{ \psfig{figure=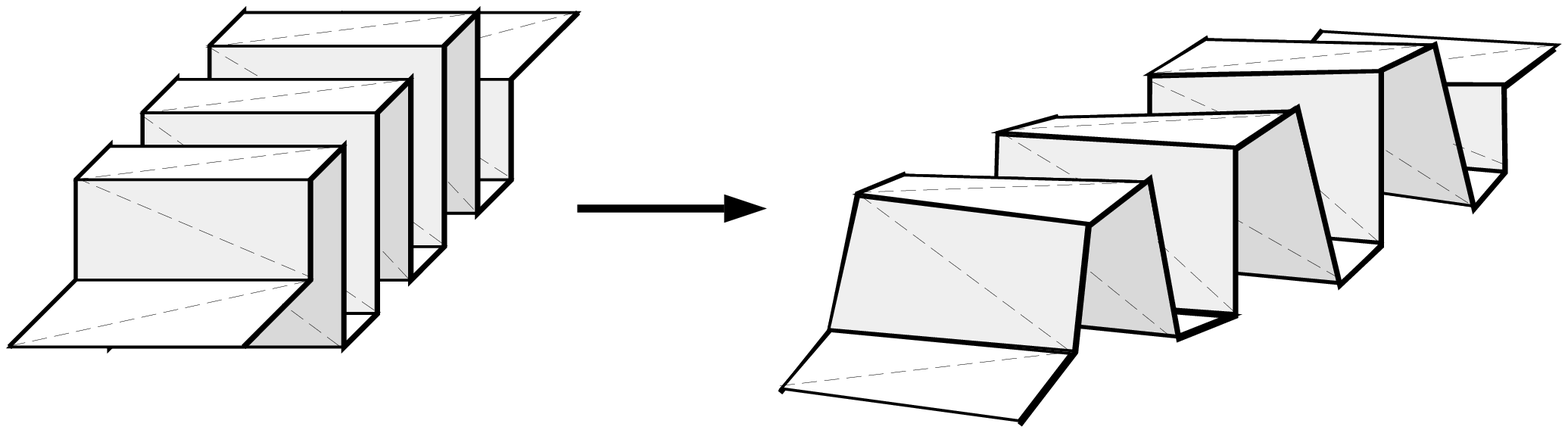,height=1.0in} }
 \mycaption{ \label{glue-edge} Interpolating map on strip.}
\end{figure}

Finally, we have to attach flexible squares to the boundary.
When we attach a flexible square to $U_2$ in Figure \ref{lin-app}
the  boundary arc lies in one of the triangular faces where 
our approximation is linear. Thus the images of the arc on the 
face of $U_2$ and  along the edge of the flexible square are 
only small distortions of each other. Thus, just as above, we 
may glue the mappings along a strip.  See Figure \ref{glue-to-edge}.

\begin{figure}[htbp]
\centerline{ \psfig{figure=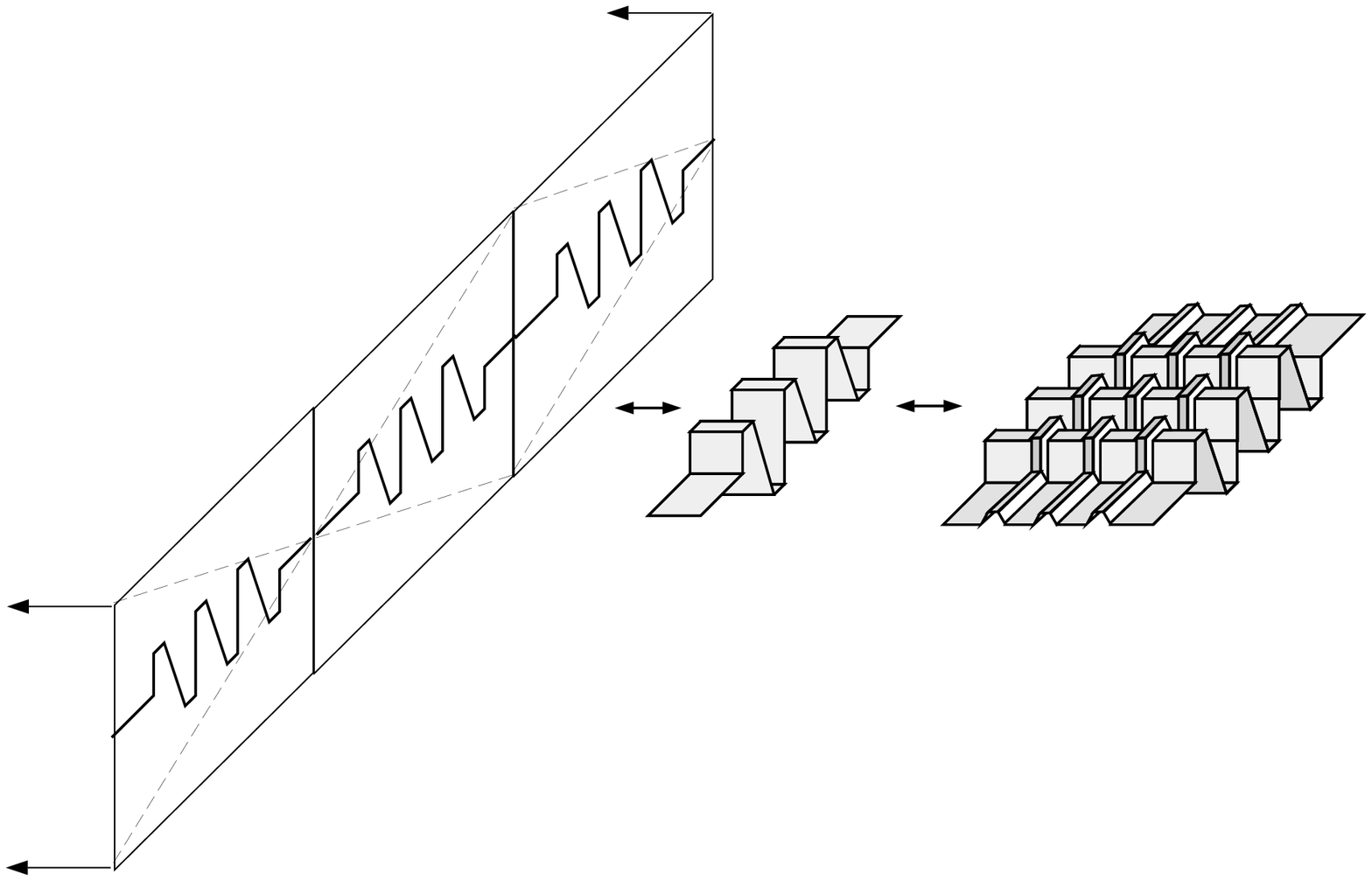,height=2.5in} }
 \mycaption{ \label{glue-to-edge} Gluing map on flexible square to 
   the boundary values.}
\end{figure}

The proof of (3) is almost exactly as above. The only real difference
is that now we may also dilate the surface by Euclidean similarities
(which change the biLipschitz constant, but not the quasiconformal
constant). Subdivide the  square into much smaller squares so  that
the Jacobian of $g$ is almost constant on each square, and replace these
squares by flexible surfaces. Then on each piece of the 
surface, we can approximate 
$g$ by the composition of a Euclidean dilation 
with the same Jacobian and a locally  biLipschitz map on the
surface. The definitions on adjacent squares can be matched as
before, so this gives the desired approximation.
\end{pf}

\begin{figure}[htbp]
\centerline{ \psfig{figure=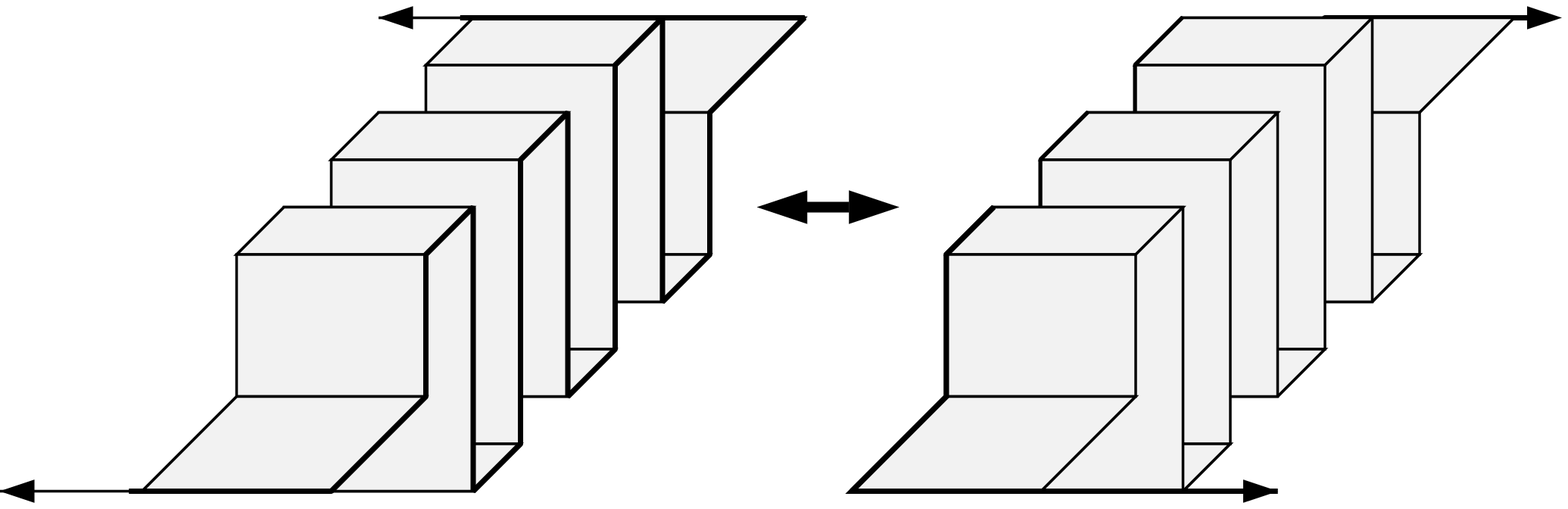,height=1.0in} }
 \mycaption{ \label{join-par} Parallel flexible squares can be 
joined.}
\end{figure}

In addition to building flexible surfaces which approximate 
a flat square, we will also want to build flexible surfaces
which approximate more complicated surfaces in ${\Bbb R}^3$.
For our purposes it will be enough to consider surfaces 
which are unions of  dyadic squares, each of which is parallel 
to one of the three coordinate planes. It is easy to 
join flexible surfaces which approximate adjoining squares
in the same plane, because the boundaries of the flexible
surfaces match exactly. See Figure \ref{join-par}.

\begin{figure}[htbp]
\centerline{ \psfig{figure=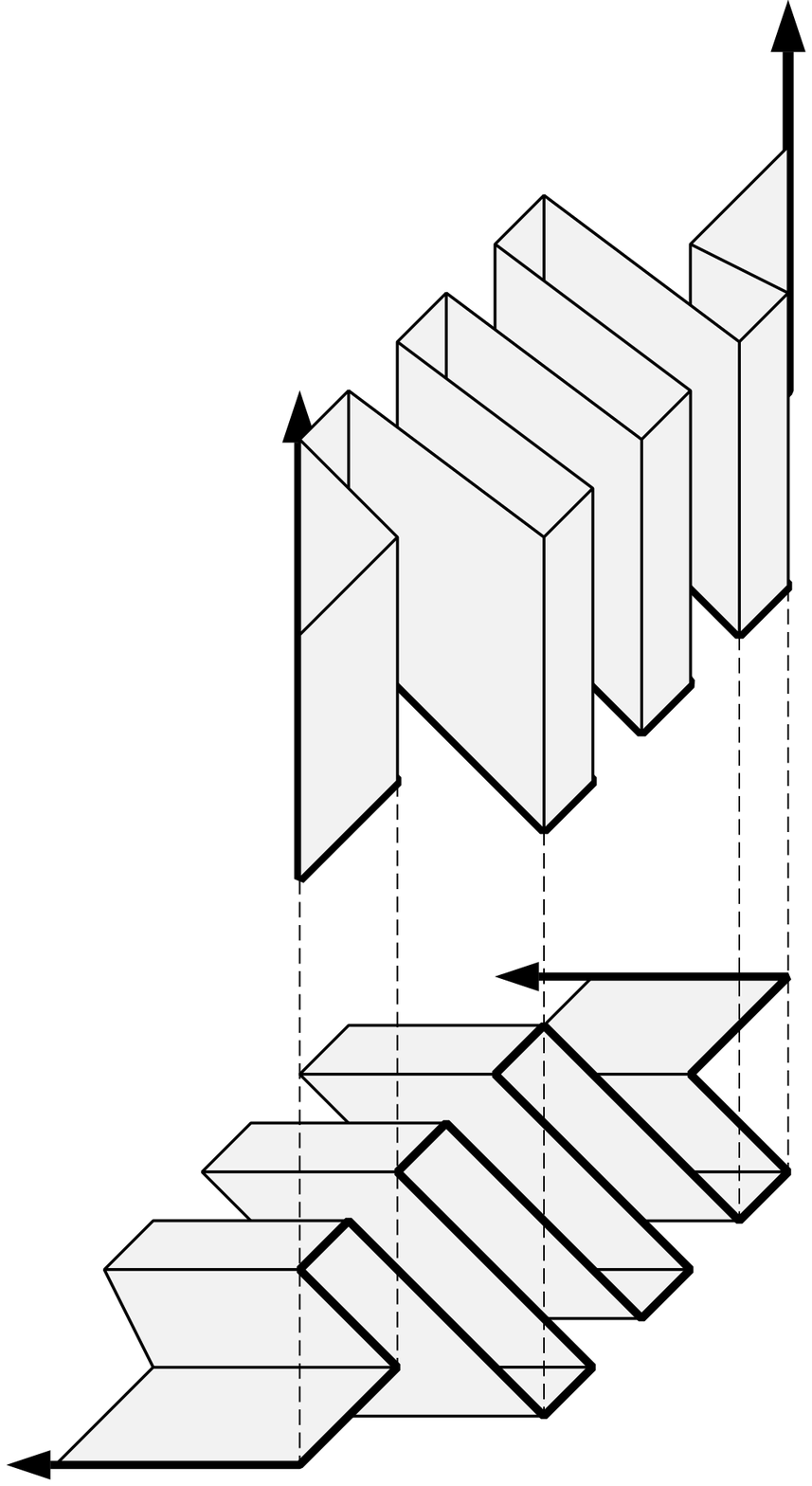,height=2.0in} }
 \mycaption{ \label{join-perp} Perpendicular flexible squares can be 
joined.}
\end{figure}

A little more care is needed if the squares are not in the same
plane. We choose three flexible arcs $\gamma_1$, $ \gamma_2 $
and $\gamma_3$ of vastly different scales,
 one corresponding to each of the coordinate
directions and we use them to  build three of flexible squares, (one
for each of the $xy$, $yz$ and $zx$ planes) with the property
that edges of these surfaces which are parallel to the
given coordinate axis have the corresponding flexible arc 
as boundary.
Thus whenever we want to join flexible surfaces corresponding 
to adjacent, but perpendicular, squares the corresponding edges 
will look the same and can be joined as in Figure \ref{join-perp}
by beveling each of the surface at 45 degrees in order to join
them.


\section{A non-removable set for locally biLipschitz maps in ${{\Bbb R}}^3$}
\label{biLip-space}

The procedure in Section \ref{biLip-plane}
can now  easily be adapted to construct a 
totally disconnected set $E$ and a homeomorphism $f$ of 
${\Bbb R}^3$ to itself which is locally biLipschitz 
off $E$, but not even quasiconformal on all of ${\Bbb R}^3$.

Let  $E^0$ be  the unit sphere in ${{\Bbb R}}^3$, let  $\Omega^0$
denote its unbounded complementary component and let $f^0$ be the 
identity map on $\Omega^0$.
Extend $f^0$  to a diffeomorphism (which we also call $f^0$)
 of ${\Bbb R}^3$ in any way you want (later we will want the 
extension to have a lot of distortion).
In general, we assume we  have a set  $E^n$ which  is a union of $J_n$ 
components $\{ E^n_j\}$, each a smooth topological  $2$-sphere
 of diameter $\leq 2^{-n}$ bounding a topological $3$-ball
$\Omega^n_j$, $j=1, \dots, J_n$.
 Denote the  unbounded complementary component of $E^n$ by  $\Omega^n_0$.
Also assume we are given a homeomorphism $f^n$ of ${\Bbb R}^3$ which 
is  locally biLipschitz on $\Omega^n_0$ and so that 
$\diam(f^n(E^n_j)) \leq 2^{-n}$.

We now describe the induction step.
 Let ${\cal S}^n$  be a collection of $\rho_n \times \rho_n$ 
cubes chosen from the usual lattice   which covers a neighborhood
 $U^n$ of $\overline{\cup_j \Omega^n_j}$. The size $\rho_n$ should be chosen 
so that $\rho_n \leq \frac 1{10} \min_{j \ne k} \dist (E^n_j, E^n_k)$, and so 
that $\diam(f^n(Q)) \leq 2^{-n}/10$.

\begin{figure}[htbp]
\centerline{ \psfig{figure=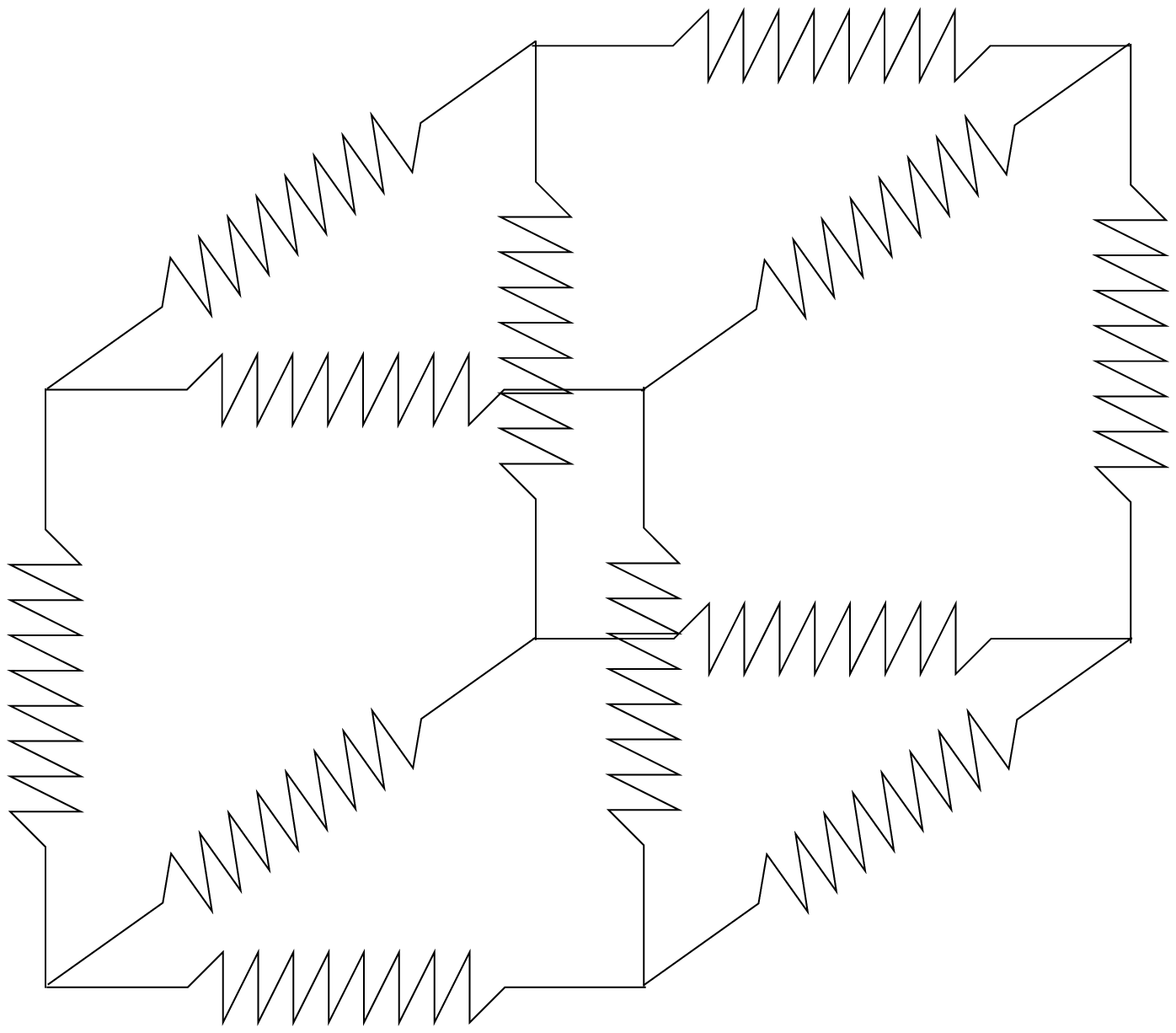,height=1.5in} }
 \mycaption{ \label{wiggly-cube} Replace edges of cubes by flexible arcs   } 
\end{figure}

Replace each  interior edge of the union of cubes   by a oscillating 
curve (such as in Lemma \ref{bad-arc-lem}) and approximate $f^{n}$
on a neighborhood of this arc by a locally biLipschitz mapping 
with a uniform constant. See Figure \ref{wiggly-cube}.

For any face of a cube in $U^n$,  take the four corners and 
consider the corresponding closed  polygonal curve obtained from the 
union of the four arcs described above. Span this curve  by a
polyhedral surface $S_1$ which approximates the original cube face and
which is a union of faces of the  dyadic squares (e.g. take a smooth spanning 
surface and replace it by faces of dyadic cubes which hit it).
Extend the map $f^{n+1}$ from the neighborhood of the 
 boundary curve to a neighborhood  of the spanning surface.
The extension should be a diffeomorphism which approximates $f^n$.
Now replace each square in the spanning surface 
by a copy of the flexible square constructed in Section 3. This 
gives a surface $S_2$.
 Then our approximation to $f^n$ on the spanning surface $S_1$ 
 has a uniformly locally biLipschitz approximation  on  a neighborhood
of $S_2$.



 Let $\Omega^{n+1}_0$ be the union of $\Omega^n \setminus 
U^n$ and the open neighborhoods constructed above on which 
$f^{n+1}$ is defined. Then $f^n$ 
extends from $\Omega^n_0\setminus U^n$ to  a locally biLipschitz
map $f^{n+1}$ on $\Omega^{n+1}_0$. Let $E^n = \partial \Omega^{n+1}_0$.
 Without loss of generality we may take $E^{n+1}$ to be a finite
union of smooth surfaces.  
Finally, extend $f^{n+1}$  to a diffeomorphism of 
of ${\Bbb R}^3$ in any way you want.


In the limit we
obtain a homeomorphism  $g$ of ${\Bbb R}^3$ to itself which is 
locally biLipschitz except on some totally disconnected set $E$.
There is enough freedom in choosing the extensions at each stage
that we can easily make sure that $f$ is not  H{\"o}lder,
so we are done.
We have now constructed a totally disconnected, non-removable set
for locally biLipschitz (and hence for quasiconformal) mappings.
The remainder of the paper deals with modifying the construction 
in order to make $E$ and $f(E)$ small in the sense of Hausdorff 
measure.

\section{Small non-removable sets for quasiconformal maps in ${\Bbb R}^2$}
\label{small-QC-2}

We now begin the process of modifying the construction so that 
$E$ and $f(E)$ are small, i.e., fix  a function $\varphi(t) = o(t^2)$, 
and show $\Haus^\varphi(E)  = \Haus^\varphi(f(E)) =0$. This 
is considerable easier for quasiconformal than for locally 
biLipschitz mappings, so we begin with a discussion of the 
quasiconformal case.  The construction in ${\Bbb R}^3$ only requires 
one extra idea, so we will first give the details in ${\Bbb R}^2$.

As before, we will define our exceptional set $E$ as a limit
of sets $\{E^n\}$, each of which is a finite union of smooth
curves $E^n = \cup_j E^n_j$
 (with diameters tending to $0$ as $ n \to \infty$).
Our  homeomorphism $f$ will be a limit of mappings 
$\{f^n\}$ which are quasiconformal on each of the finitely 
many components  of $ \Omega^n = {\Bbb R}^2 \setminus E^n$.
This is different from what we did before, where we only defined
$f^n$ to be ``good''  on the single unbounded component and 
extended it any we wanted to the bounded components.
The maps $\{ f_n \}$ will  not be  homeomorphisms because 
the definitions on different components of $\Omega^n$
 will disagree on $\partial \Omega^n=E^n$.
The main idea of the inductive step is to 
reduce the amount of disagreement at each step.

\begin{figure}[htbp]
\centerline{ \psfig{figure=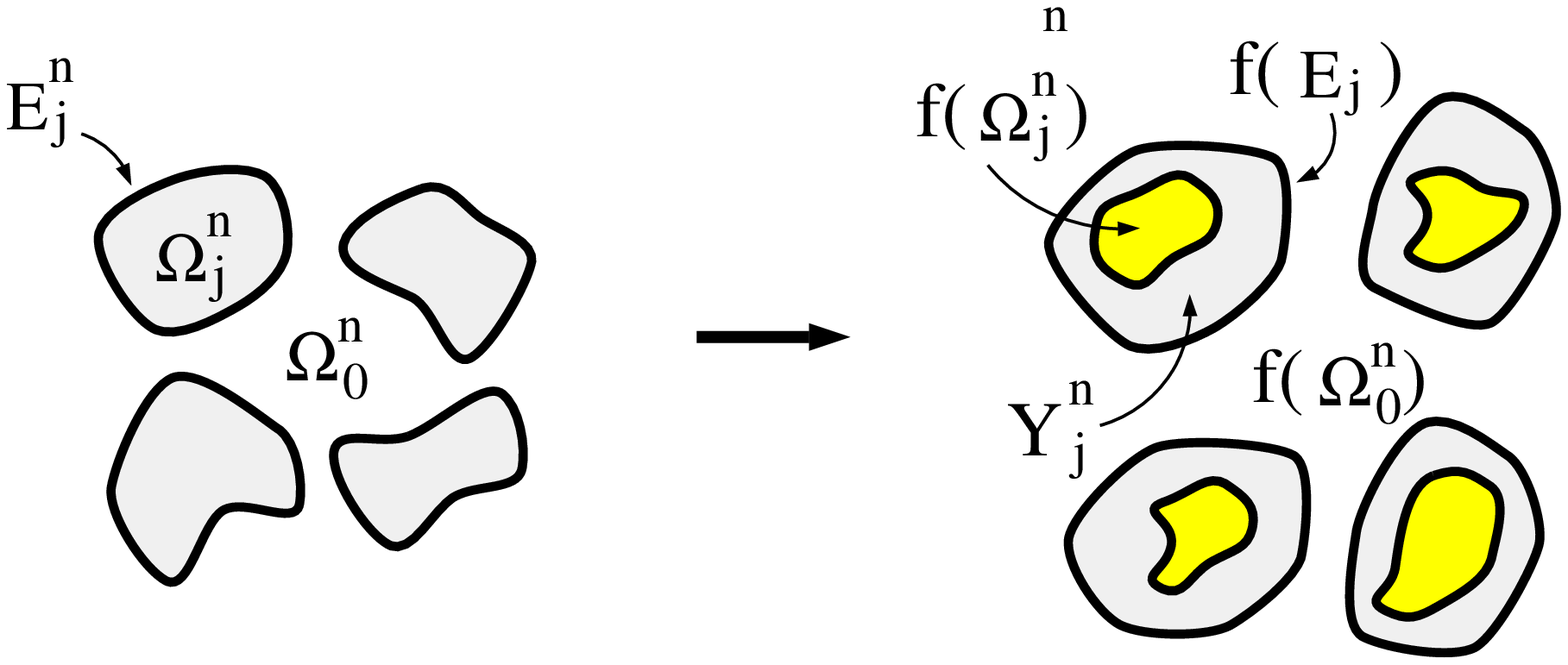,height=1.5in} }
 \mycaption{ \label{ind-hyp-def}  Definitions at the beginning of 
the induction.}
\end{figure}

To begin the induction we start with $E^0 = \{z: |z|=1\}$.
 Let $\Omega^0$ be the complement of $E^0$,   let $\Omega^0_0$ 
be its unbounded component and $ \Omega^0_1$ the 
bounded component.
Define 
$$ f^0_0(z) = z, \quad z \in \Omega^0_0,$$
 $$ f^0_1(z)= \frac 12 z, \quad  z \in \Omega^0_1.$$
The induction hypothesis is as follows.
Suppose we are given a compact set $E^n $
 which is a finite
union of $ J_n$ smooth closed curves, $\{ E^n_j\}$, which are 
disjoint with disjoint interiors. 
Let $\Omega^n$ be the complement
of $E^n$. Its unbounded component is denoted $\Omega^n_0$ and 
the bounded components are denoted $\Omega^n_j$, $j=1,\dots,J_n$.
Suppose we are given diffeomorphisms $f^n_j$ on $\Omega^n_j$, $j=0,
\dots, J_n$, which are
quasiconformal with constant $M$ on each component, and  
$f^n_j(\Omega^n_j)$ lies 
in  $Y^n_j$, the bounded complementary component of $f^n_0(E^n_j)$.
See Figure \ref{ind-hyp-def}

What follows is  a description of how to construct $f^{n+1}$
and $E^{n+1}$ from $f^n$ and $E^{n}$.

\begin{figure}[htbp]
\centerline{ \psfig{figure=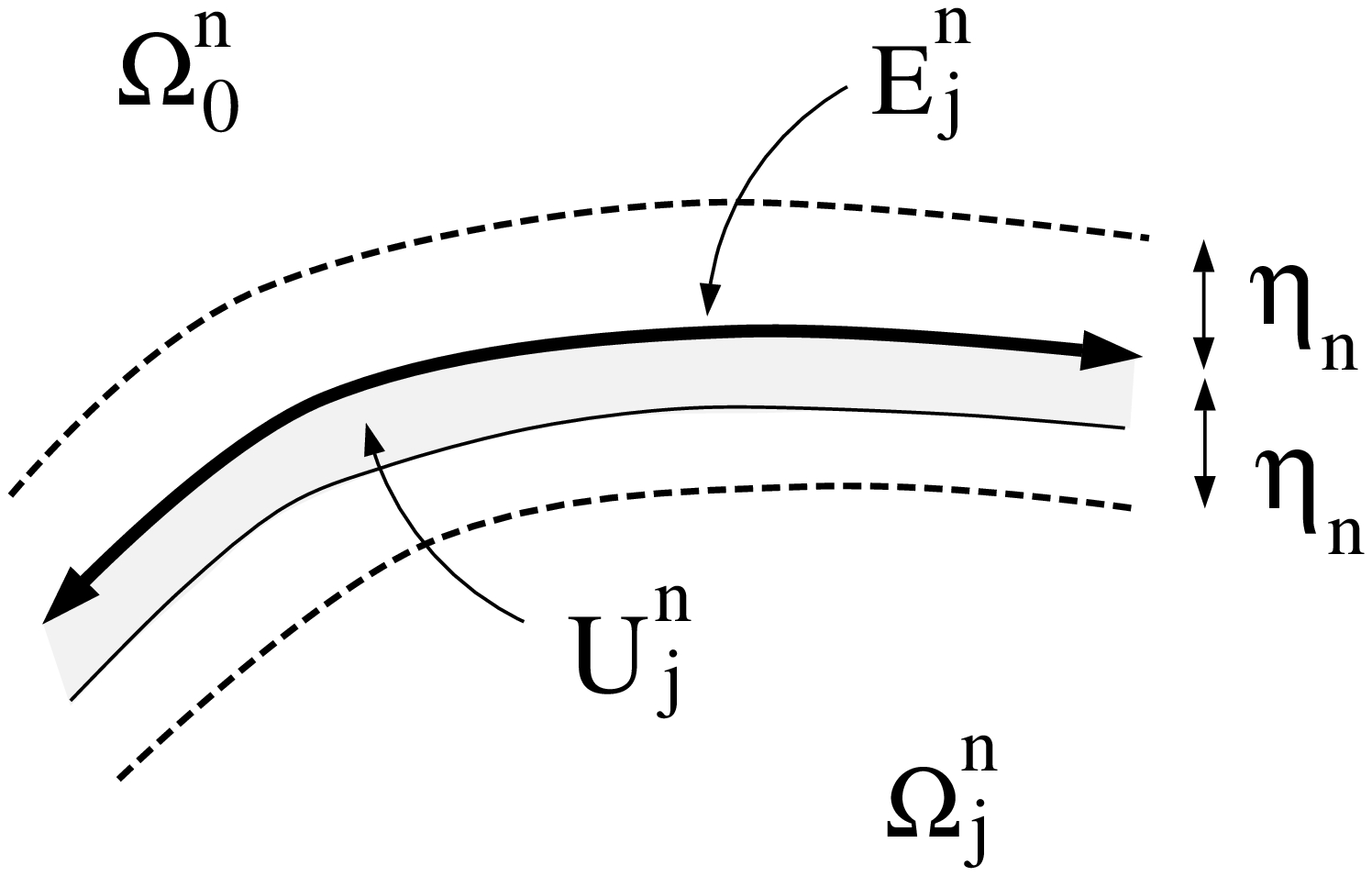,height=1.5in} }
 \mycaption{ \label{def-U}  Definition of $U^n_j$.}
\end{figure}

\begin{figure}[htbp]
\centerline{ \psfig{figure=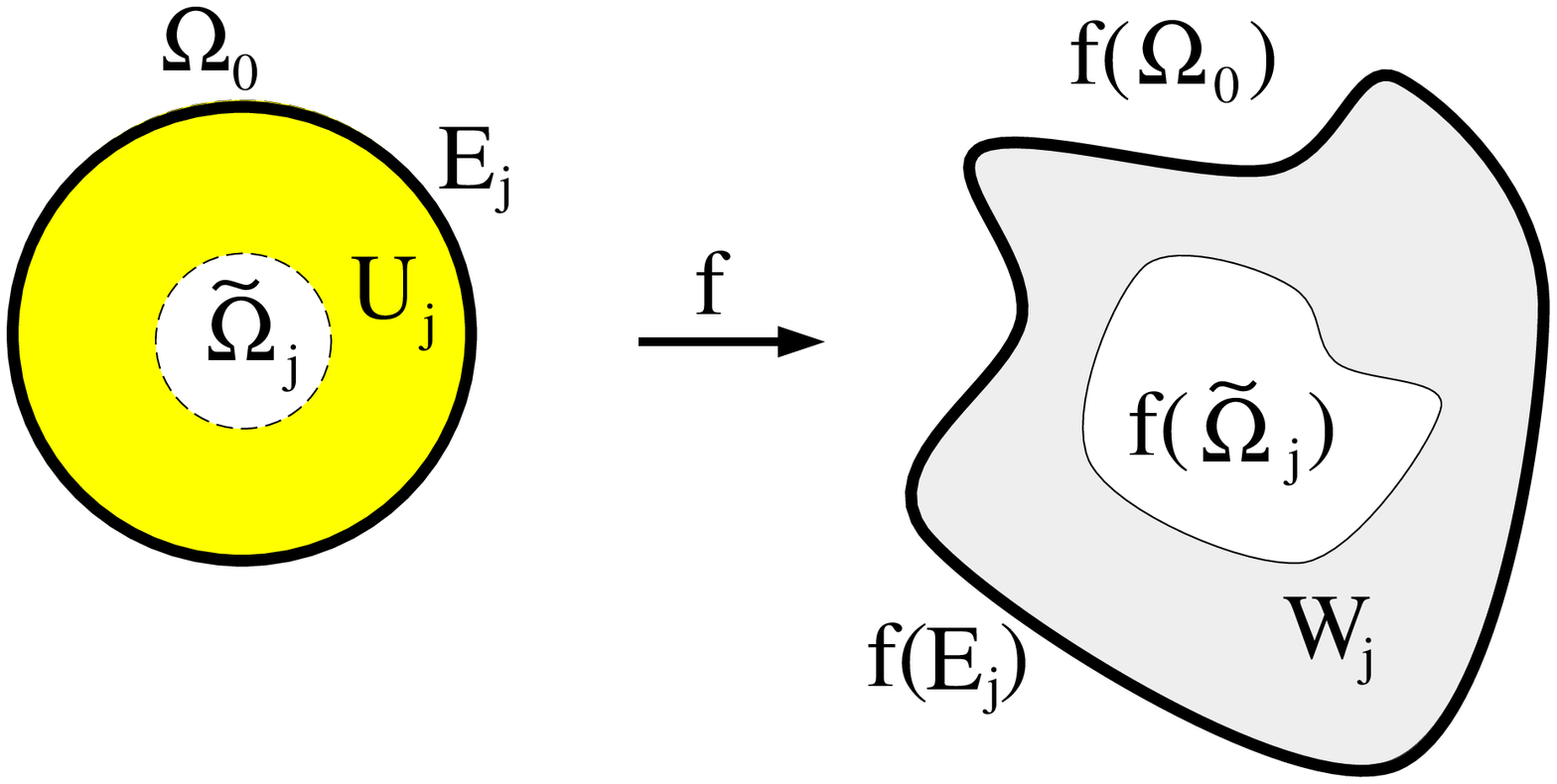,height=1.5in} }
 \mycaption{ \label{defns}  Definitions of $U^n_j$ and $W^n_j$.}
\end{figure}

\begin{figure}[htbp]
\centerline{ \psfig{figure=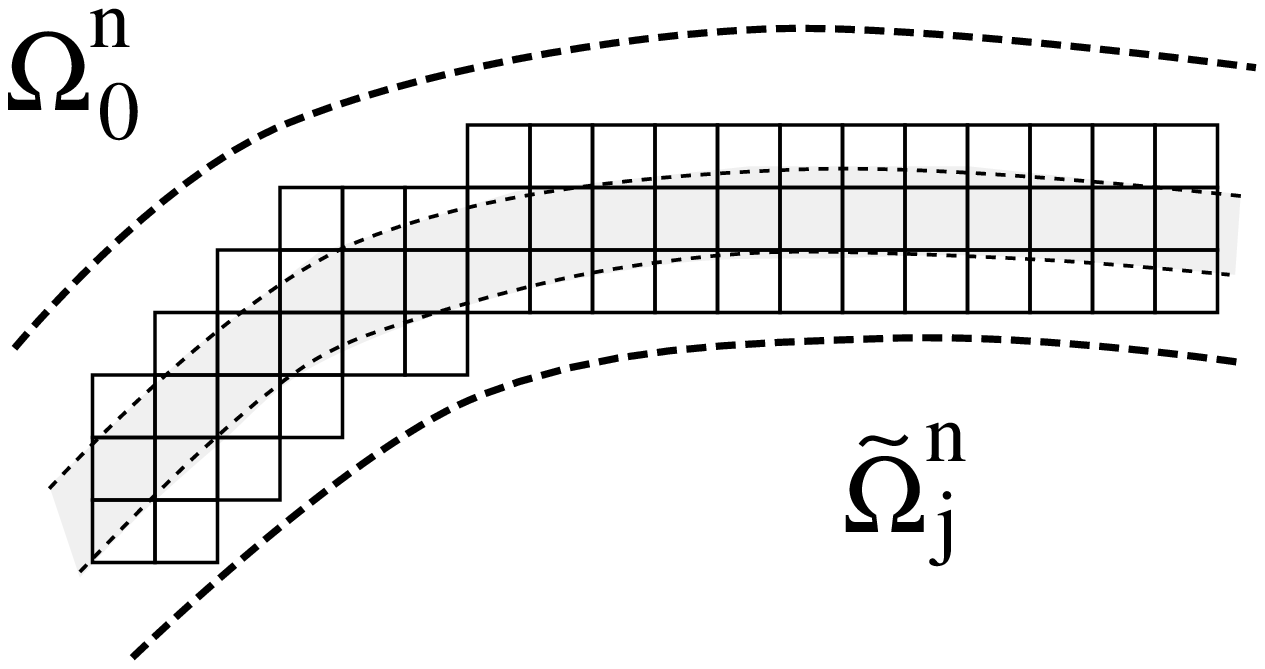,height=1.20in} }
 \mycaption{ \label{def-Sn}  Definition of ${\cal S}_n$.}
\end{figure}

{\bf Step 1:}
Fix a  very small number $\eta_n>0$ 
and  for $j=1, \dots, J_n$ let
$V^n_j$ be a smooth  annular neighborhood of $E^n_j$ which is 
contained in  
$$  \{ z :  \dist(z, E^n_j) < \eta_n/10\}
  \subset \widehat U^n_j  \equiv  \{ z :  \dist(z, E^n_j) < \eta_n\},$$
and let $U^n_j = V^n_j \cap \Omega^n_j$.
Let $U^n =\cup_j U^n_j$.
See Figure \ref{def-U}.
 For $j=1, \dots, J_n$, let 
$$ \widetilde \Omega^n_j = \Omega^n_j \setminus \overline{U^n}, \qquad
 W^n_j  =  Y^n_j \setminus \overline{ f^n_j (\widetilde \Omega^n_j)) }.$$
See Figure \ref{defns}.
Then $W^n = \cup_j W^n_j$ consists of $J_n$ annuli, so there is smooth 
diffeomorphism from  $U^n$ to $W^n$, which maps $U^n_j $ diffeomorphically
to $W^n_j$ and which agrees with $f^n_0$ on $E^n$ and with 
$f^n_j$ on $E^n_j$.
Thus we can construct a smooth diffeomorphism 
$g^n: {\Bbb R}^2 \to {\Bbb R}^2$ which agrees  with $f^n$ on 
${\Bbb R}^2 \setminus U^n$.
 Now choose $\rho_n < \eta_n/100$ and consider the 
grid of $\rho_n \times \rho_n$ squares  from the 
lattice $\rho_n {\Bbb Z} \times  \rho_n{\Bbb Z}$.
Let ${\cal S}_n$ be  a   collection of such squares which cover
$U^n$ and  are  contained in  $\widehat U^n$.
Let $F= \cup_{Q \in {\cal S}} \partial S$.  
See Figure \ref{def-Sn}.

{\bf Step 2:}
We  want to extend $f^n$ from $\Omega^n$ to an open connected set 
$\Omega^{n+1}_0$ which contains $\Omega^{n}_0$, $\cup_j 
\widetilde \Omega^{n}_j$ and a neighborhood of $F$ by 
approximating $g^n$ on a neighborhood of $F$.
 We can do this because any mapping on
a line segment can be approximated by a mapping with bounded 
quasiconformal distortion on a neighborhood of the interval.
In our case, it 
is very simple to draw a picture of the approximations.
(We use straight lines instead of flexible arcs so that the
resulting domain will be quasiconvex. See Section \ref{quasi-con}.)

In a neighborhood of a  corner $x$  of $F$ we simply define 
$f^{n+1}$ to be a Euclidean similarity with the property
 that $f^{n+1}(x) = g^n(x)$. 

\begin{figure}[htbp]
\centerline{ \psfig{figure=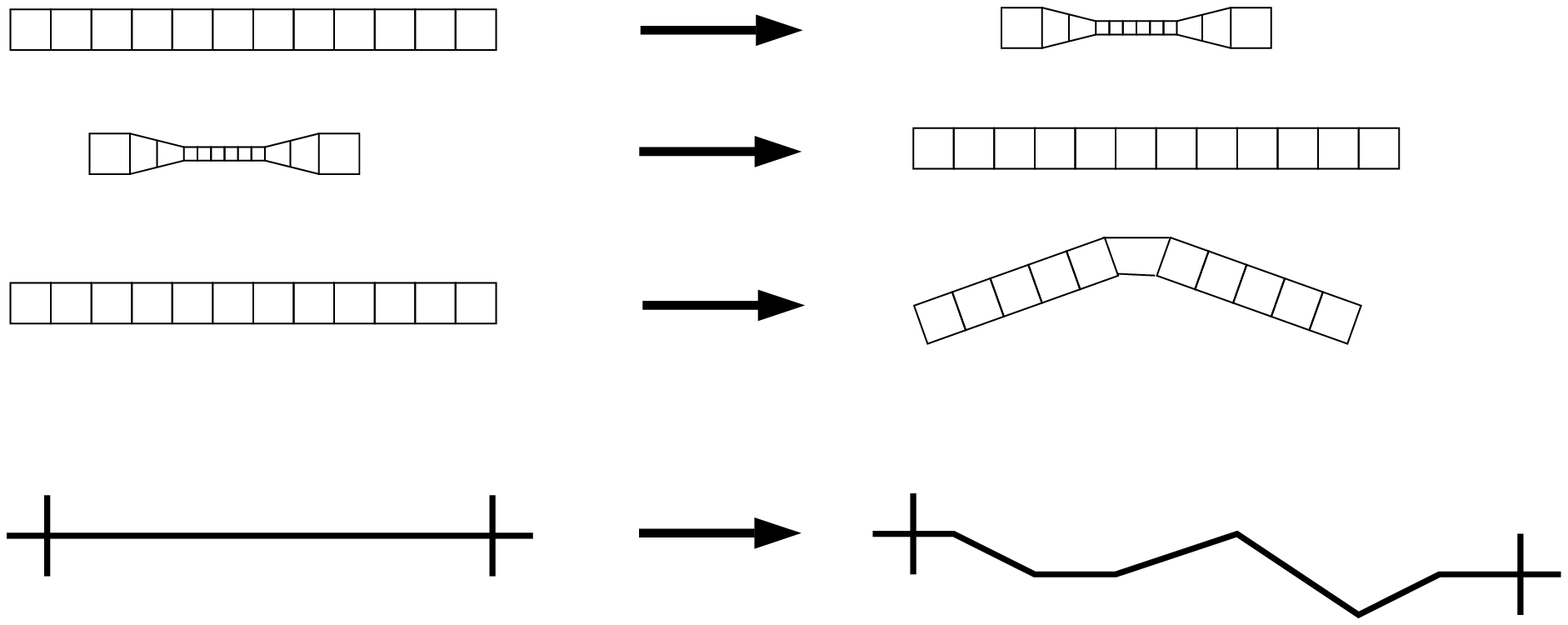,height=2.0in} }
 \mycaption{ \label{bowtie}   
Stretching or shrinking a line segment to approximate a polygonal arc }
\end{figure}

On the  line segments connecting corners we approximate 
$g^n$ by a quasiconformal map. Figure \ref{bowtie} shows how line segments 
may be stretched, shrunk or bent by means of a quasiconformal 
map with uniformly bounded dilation.
Thus by approximating $g^n$ by a polygonal arc and 
using these maps to approximate each segment, we 
obtain the desired map.

\begin{figure}[htbp]
\centerline{ \psfig{figure=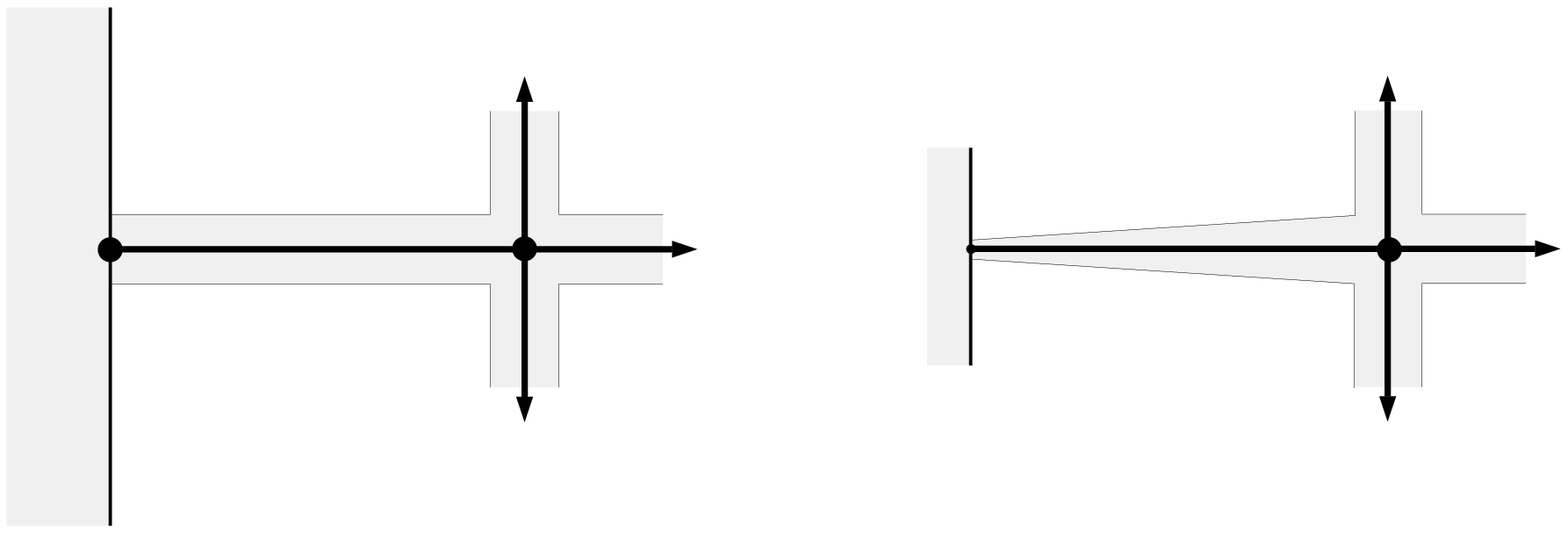,height=1.5in} }
 \mycaption{ \label{connect} Make the connections between $f^{n+1}$ and 
$f^n$} 
\end{figure}

The only remaining observation we have to make is that the 
approximation can be chosen to agree with the map $f^n$
outside $U^n$.
Suppose $x$ is a ``boundary corner'' of $F$.
Then $x$ is 
connected by one or more grid segments  to points in 
$ \widetilde \Omega^n$. On a  segment connecting $x$ and $y$
 we define $f^{n+1}$ to map the arc 
so that $f^{n+1}$  extends 
both $f^n$ and $f^{n+1}$  to a neighborhood of  the line segment. 
See Figure \ref{connect}.

\begin{figure}[htbp]
\centerline{ \psfig{figure=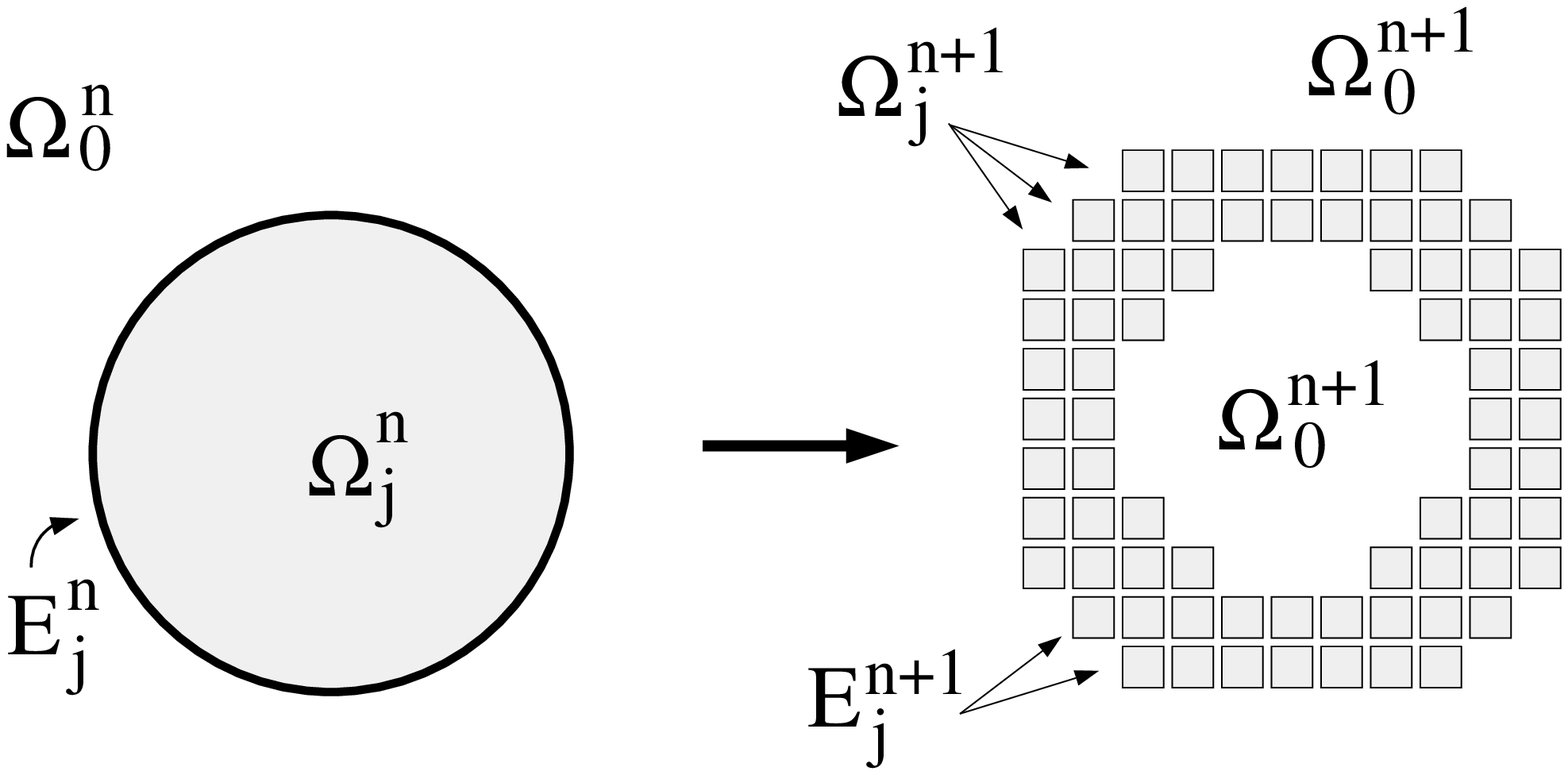,height=2.0in} }
 \mycaption{ \label{replace-pic} $E^n$ and $E^{n+1}$ } 
\end{figure}

We now have a smooth diffeomorphism $f^{n+1}$ defined on an open set
$\Omega^{n+1}_0$ which contains $F \cup \Omega^n_0  \cup_{j\geq 1} 
 \widetilde \Omega^n_{j}$. 
Without loss of generality we may assume that $\partial \Omega^{n+1}_0$ is 
a finite union of smooth closed curves. Let 
$E^{n+1} = \partial \Omega^{n+1}_0$, and let $\Omega^{n+1}_j $, $j >1$ be an 
enumeration of the finitely many bounded complementary 
components. 
See Figure \ref{replace-pic}.
To avoid confusion, let $f^{n+1}_0$ denote the 
continuous extension of $f^{n+1}$ from $\Omega^{n+1}_0$ to its 
closure.

\begin{figure}[htbp]
\centerline{ \psfig{figure=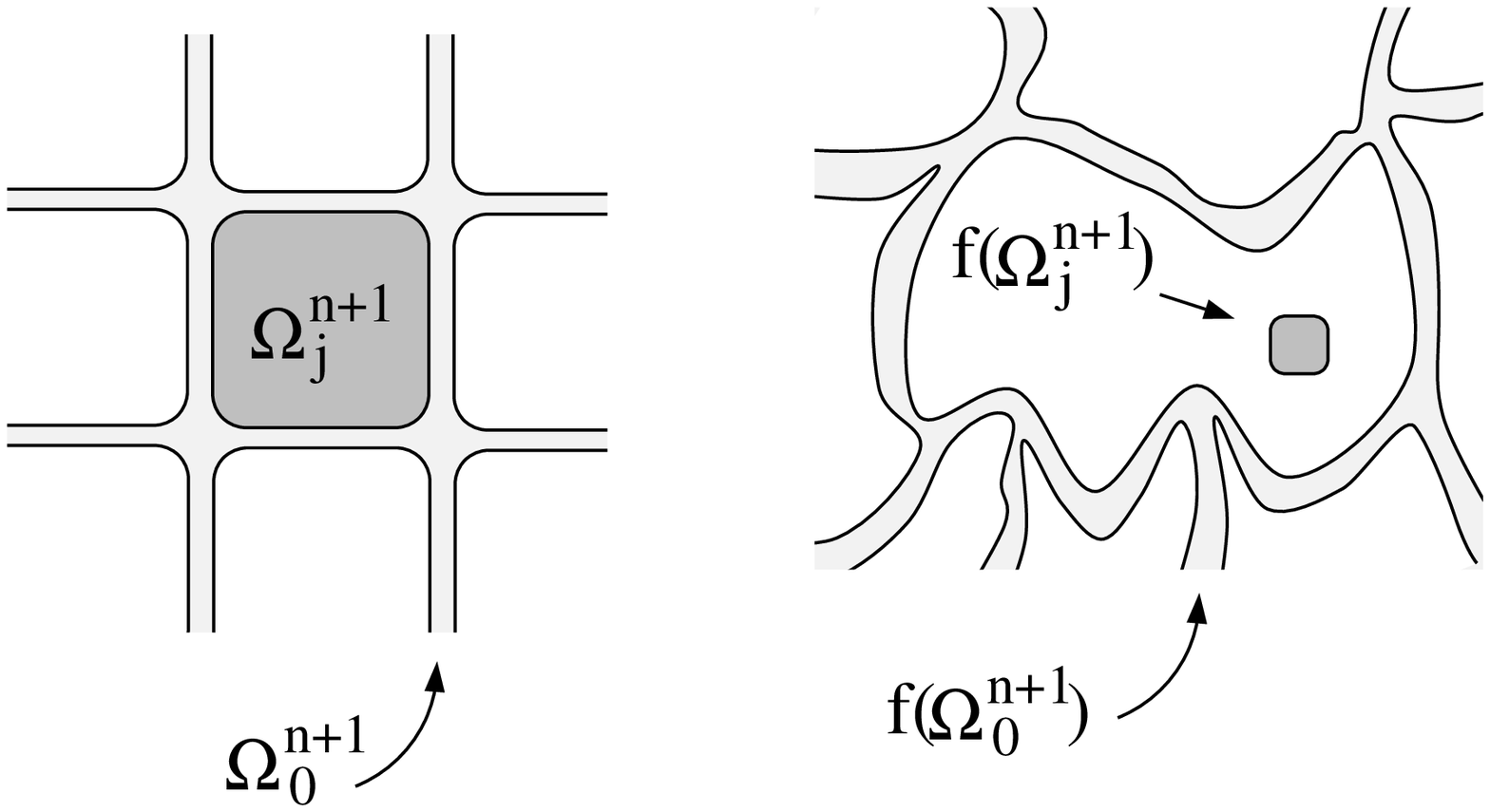,height=2.0in} }
 \mycaption{ \label{sim-map} Defining $f^{n+1}$ on the bounded components } 
\end{figure}

 To define $f^{n+1}_j$ on $\Omega^{n+1}_j$ we simply choose it 
to be a Euclidean (orientation preserving) similarity which 
maps $\overline{\Omega^{n+1}_j}$ into $Y^{n+1}_j$, the region bounded by 
$f^{n+1}_0(E^{n+1}_j)$. 
See Figure \ref{sim-map}.


This completes the inductive step of the construction, i.e. given the
 set $E^n$ and mapping $\{f^n_j\}$ we have constructed $E^{n+1}$ and 
$\{ f^{n+1}_n\}$ which satisfy the induction hypothesis.
In particular,  if we let $F^n = \overline{\widehat U^n}$
and let $\{ g^n\}$ be the maps constructed at the end of 
Step 1, they satisfy Lemma \ref{nested-sets}.
Using the lemma we see our maps converge to a homeomorphism
which is clearly quasiconformal off a Cantor set $E$.
Finally, to see that $g$ is not quasiconformal on all of ${\Bbb R}^2$, 
there are several things we could do. The easiest is to define the
homeomorphisms $\{g^n\}$ in Step 1 so that the limiting 
homeomorphism  $f$ is not H{\"o}lder of any positive order.

To see that $E$ can be taken to have $\Haus^\varphi(E)=0$, fix 
a function $\varphi(t) = o(t)$. Since $E^n$ is a finite 
union of smooth curves, it has finite length and can be covered 
by $C_n r^{-1}$ disks $\{ D_j\}$  of size $r$ (for all small enough $r$).
Choose $r_n$ so small that $\varphi(r_n) <  \frac 1n C_n^{-1} r_n$.
In Step 1 of the construction choose $\eta_n << r_n$ so small that 
$\widehat U_n \subset \cup_j D_j$. From the construction it is clear
that we can take $E \subset  \cup 2D_j$, so  
$$ \Haus^\varphi(E) \leq \lim_{n\to \infty}  2 \sum_j \varphi(r_n)
 \leq \lim_{n\to \infty}   \frac 2n =0,$$
as desired.

 We can define
the neighborhoods of $F$ to be so small at each stage that 
$\text{area}(W_n)$ remains bounded away from $0$ for all $n$.
This means that $f(E)$ can have  positive area.

If we want to make $f(E)$ small, then instead of 
defining $f^{n}$ to be a similarity on $ \Omega^{n}_j$, 
define it to be a conformal mapping from $\Omega^{n}_j$ (which 
is topologically a disk) to $Y^{n+1}_j$ (which is also a disk).
Then at the next step the annular regions $W^{n+1}_j$ can be taken 
to lie in an arbitrarily thin  neighborhood of $f^{n+1}_0(E^{n+1})$.
By taking a small enough neighborhood we can obtain 
$\Haus^\varphi(f(E)) =0$, just as above.
This last step (where we have used the Riemann mapping theorem) 
is the only  one which causes a problem in ${{\Bbb R}}^3$.

\section{The quasiconformal construction in ${\Bbb R}^3$}
\label{small-QC-space}

 As before let $E^0$  
be the unit sphere, $\Omega^0 = \Omega^0_0 \cup \Omega^0_1$ its
complement and  
$$ f^0_0 (z) =  z, \quad z \in \Omega^0_0,$$
$$ f^0_1 (z) =  \frac 12 z, \quad  z \in \Omega^0_1.$$

In general, suppose we have a compact set $E^n$ consisting 
of $J_n$ components $\{ E^n_j\}$, each of which is a smooth surface 
diffeomorphic to  the 2-sphere and bounding a topological 
3-ball $\Omega^n_j$. Let $\Omega^n_0$ be the unbounded complementary
component of $E^n$. Assume we have a quasiconformal map $f^n_j$, 
$j =0 , \dots, J_n$ defined on each component. These maps 
 extend smoothly across the boundaries and 
$f^n_j(\Omega^n_j) $ is a subset of $Y^n_j$,  the bounded  complementary 
component of $f^{n}_0(E_j^n)$.

{\bf Step 1:}  Define an open set   $U^n_j \subset \Omega^n_j$
which is a topological annulus
 (i.e., homeomorphic to 
$S^2 \times (0,1)$) with one boundary component $E^n_j$ and so 
that $U^n_j$ lies in 
a $\eta_n/10$ neighborhood of $E^n_j$. Let $U^n = \cup_j U^n_j$.
 For $j = 1, \dots, J_n$, let 
$$ \widetilde \Omega^n_j = \Omega^n_j \setminus U^n , \qquad 
 W^n_j = Y^n_j \setminus f^n_j(\widetilde \Omega^n_j).$$
Then $W^n_j$ is an  annulus  (see  Remark 6.1
concerning the annulus conjecture) and hence is 
diffeomorphic to $U^n_j$.
Therefore there is a diffeomorphism $g^n$ of ${\Bbb R}^3 $ to itself
 which agrees with $f^n_j$ on $ \widetilde \Omega^n_j$. 
As before choose $\delta_n < \eta_n$ and consider  a collection 
of cubes from a  $\delta_n$-grid  which covers
$U^n$ and lies in a $\eta_n/2$ neighborhood of $E_n$.
 Let $F$ denote the union of the faces of these cubes.

{\bf Step 2:}
We want to  define an approximation $f^{n+1}$ to $g^n$ on a neighborhood of 
$F$, but may be impossible. Instead we 
will define the approximation at the corners and along the edges of 
the cubes and then replace the faces by copies of our ``flexible
squares''. We then define the approximation on a neighborhood of 
these surfaces using Lemma \ref{flex-square}.

On a neighborhood of each corner we define our approximating map  $f^{n+1}$
to be a similarity which  agrees with $g^n$ at the corner point. 

On each line segment connecting two such corner points we 
define a uniformly quasiconformal approximation $f^{n+1}$  on some 
neighborhood.
This is exactly the same 
as the two dimensional case, since we can stretch or contract
a line segment by replacing the squares by cubes in Figure 
\ref{stretch-arc}.

 For each face of each cube, extend the approximation defined 
above in a neighborhood of the four edges to an approximating 
diffeomorphism defined on a neighborhood of the face (do this 
in any smooth way without worrying about the quasiconformal 
constant).

Next, the   face of each cube is replaced by a scaled copy 
of a ``flexible square''.
The surfaces are chosen so that they lie in the 
neighborhoods of the faces described in the previous paragraph. 
By construction we have a uniformly  quasiconformal approximation to 
$g^n$ on some neighborhood of these surfaces.

Let $\Omega^{n+1}_0$  be the  open set where $f^{n+1}$ has been defined
and let $E^{n+1} = \cup_j E^{n+1}_j$ be its boundary components.
Without loss of generality we may assume these are smooth. Let
$\Omega^{n+1}_j$ be the bounded complementary component of
$E^{n+1}_j$.
In each  $\Omega^{n+1}_j$,
we define a new mapping 
by a Euclidean similarity, so that the image of the component
is contained in $Y^{n+1}_j$,  the bounded component 
 of the complement of $f^{n+1}(\partial \Omega^{n+1}_j)$.

This completes the induction step. The process of passing to the limit 
is exactly as in the two dimensional case. Similarly, the proof that
$\Haus^\varphi(E) =0$ is unchanged, except that $E$ now lies
in a thin neighborhood of a surface instead of a curve, so we 
get an estimate for $\varphi(t) = o(t^2)$ instead of $o(t)$.

If we want to make $\Haus^\varphi(f(E)) =0$, the argument used 
in the two dimensional case does not work here. In that case
we defined $f^{n}$ on the bounded components to be a conformal 
mapping using the Riemann mapping theorem, but in ${\Bbb R}^3$, 
this is not available to us. However, we can achieve the same 
result by using the following observation.

\begin{lemma} \label{qc-sub}
Suppose $\Omega $ is an open connected set with a smooth 
boundary and suppose $Q= [0,1]^3$ is the unit cube. 
Then there is a quasiconformal map $h$  of $Q$ onto 
a subdomain $\widetilde \Omega \subset \Omega$ such that 
$E = \overline{ \Omega} \setminus \widetilde \Omega$  
has $\sigma$-finite 2-dimensional measure.
\end{lemma}
 
\begin{figure}[htbp]
\centerline{ \psfig{figure=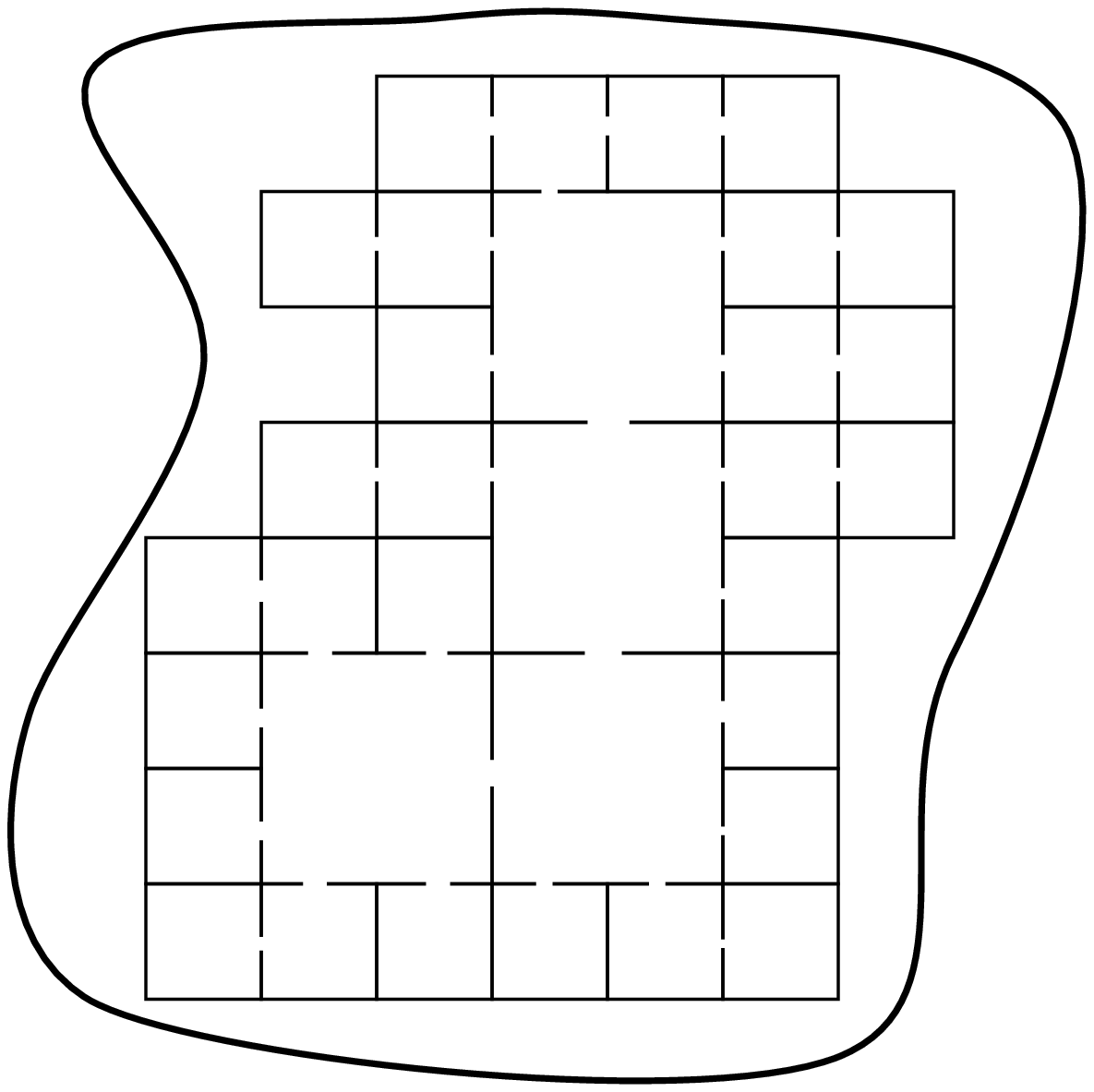,height=1.75in} }
 \mycaption{ \label{whit-union}   A union of Whitney cubes with ``openings''}
\end{figure}

\begin{figure}[htbp]
\centerline{ \psfig{figure=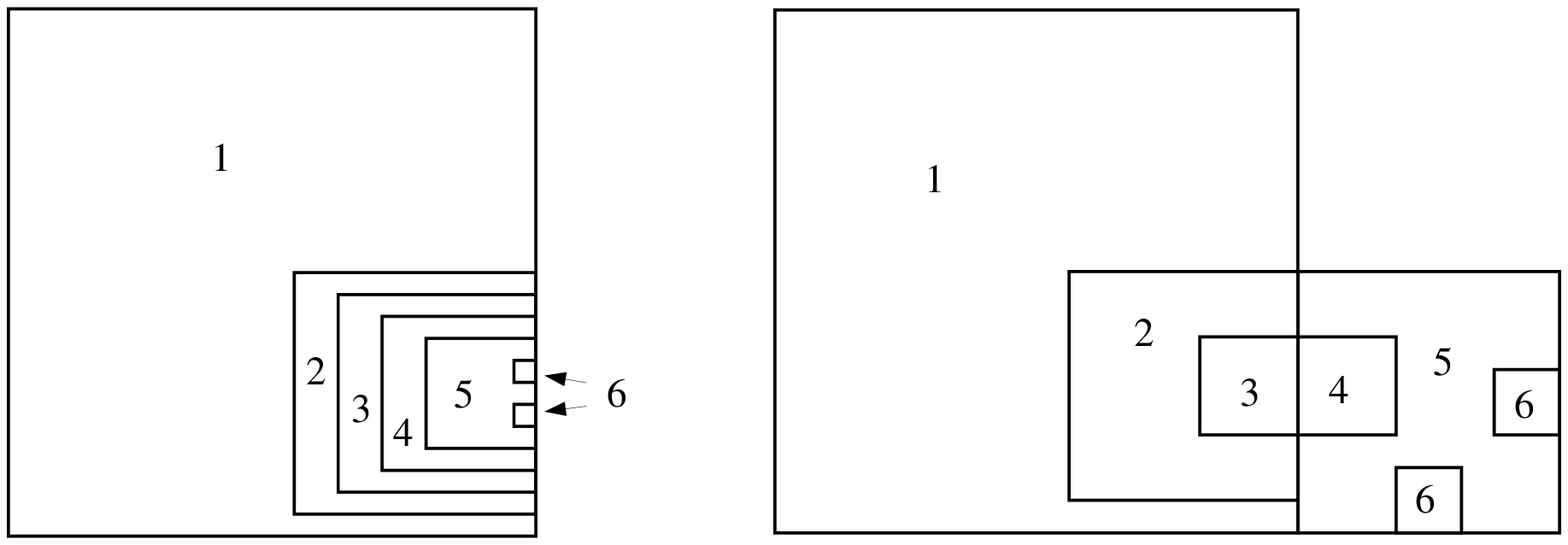,height=1.75in} }
 \mycaption{ \label{QCmap-nums} Mapping a cube to a union of cubes} 
\end{figure}

\begin{pf}
To prove this one simply takes a Whitney decomposition
$\{ Q_j \}$ for $\Omega$. Let $\widetilde \Omega$ be the 
union of the interiors of these cubes, plus small openings
between certain adjacent cubes. This can be done so that
$\widetilde \Omega$ is  connected and  simply connected.
See Figure \ref{whit-union}.
 It is not hard to see that $\widetilde \Omega$ is quasiconformally
equivalent to $Q$. For example, Figure \ref{QCmap-nums} shows how to map
one cube  quasiconformally to the union of two; and in such a 
way that the map is conformal where additional cubes might be 
attached.
Since $\Omega \setminus 
\widetilde \Omega$ is contained in a countable number of flat 
squares, the final claim is obvious.
\end{pf}

Using this one can get $\Haus^\varphi(f(E))=0$ (for some 
$\varphi(t) = o(t^2)$) as follows.
Instead of using a Euclidean similarity to map each 
 component $\Omega^{n+1}_j$ into the appropriate 
component, use the previous lemma applied to
$\Omega = Y^{n+1}_j$ and  
cube  $Q$ containing $\Omega^{n+1}_j$.
Fix a sufficiently small $r$ and choose
a covering of $\Omega \setminus \widetilde \Omega$ 
with $ (n \varphi(r))^{-1}$ cubes of size $r$.
If the flexible surfaces making up the faces of $\Omega^{n+1}_j$ are
close enough to the faces of $Q$, then $h(Q\setminus
\Omega^{n+1}_j)$ will be contained in a small neighborhood 
of $\Omega \setminus \widetilde  \Omega$ and will also be covered 
by these cubes. Thus we can cover all of 
$f(E) $ by only $\frac 1n (\varphi(r) r^2)^{-1}$ cubes of 
size $r$ which is enough to give $\Haus^\varphi(f(E)) =0$.



\begin{remark}
We now address the topological problem alluded to in the construction.
It concerns the statement that each $W_j^n$ is a topological
annulus.
We would like to know that given 
 a closed  $n$-ball $B \subset {\Bbb R}^n$ and a homeomorphism $h: {\Bbb R}^n \to 
{\Bbb R}^n$ with $h(B) \subset \text{int(B)}$ then $B\setminus h(B)$
is homeomorphic to $S^{n-1}\times (0,1)$.  This  may seem obvious, 
but it is known as the annulus 
conjecture and was only proven for $n=3$ by Moise in 1952 \cite{Moise}
(for $n=4$ it was proven  by  Quinn in 1982 \cite{Quinn}
and for $n> 4$ by Kirby in 1969 \cite{Kirby}).
Fortunately, in our case the $3$-balls in question are very 
explicit polyhedron and the existence of the desired homeomorphism 
is fairly clear. Moreover, our case fits into either the quasiconformal 
or biLipschitz categories and these cases are handled by work of 
Sullivan and of Tukia and V{\"a}is{\"a}l{\"a} \cite{Tukia-Vaisala}.
\end{remark}

%

\section{Quasiconvexity and product sets}
\label{quasi-con}

\begin{figure}[htbp]
\centerline{ \psfig{figure=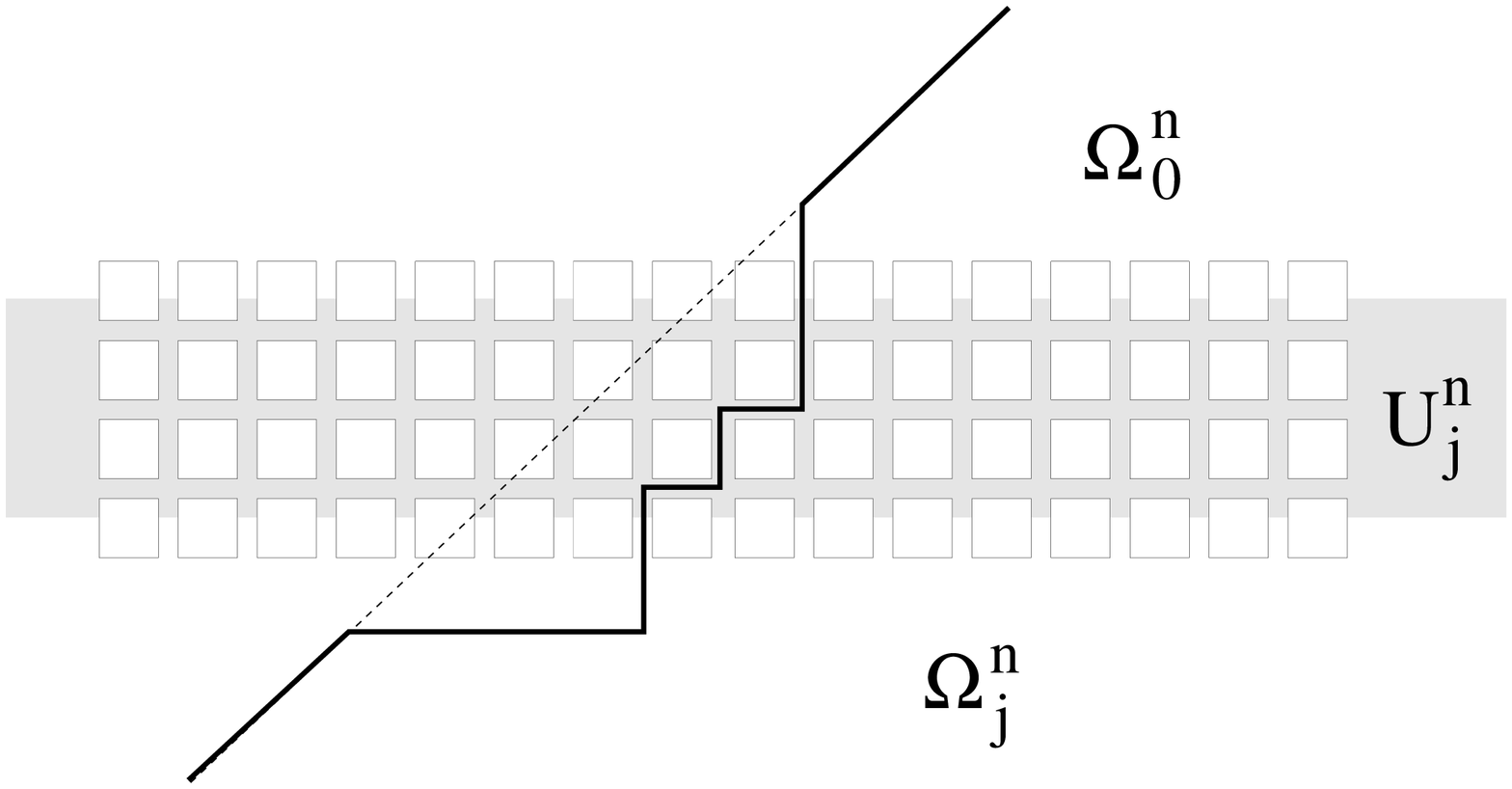,height=1.75in} }
 \mycaption{ \label{mod-arc} Proving $E$ has quasiconvex complement } 
\end{figure}

The non-removable sets for quasiconformal mappings 
constructed  in the two previous  sections are  removable for locally 
biLipschitz mappings.
To see why, we first claim that the complement 
$\Omega = {\Bbb R}^2 \setminus E$ is {\it quasiconvex}
\cite{Heinonen-Koskela96}, i.e.,  that any 
two points $x,y$ in $\Omega = {\Bbb R}^2 \setminus E$ can be connected
by a path in $\Omega$ with length $\leq C |x-y|$. 
We may assume $x$ and $y$ are both in $\Omega^{n}_0$ for 
some $n$ and  simply take the line segment between $x$ and $y$, 
except that whenever the segment crosses $E^n$ between
 $\Omega^n_0$ and one of the components $\Omega^n_j$, 
we modify it to be a polygon 
arc whose sides lie along the edges of cubes covering 
$U^n_j$. See Figure \ref{mod-arc}.
  These edges are in $\Omega^{n+1}_0 \subset 
\Omega$ by construction and the modification at most 
doubles the length of the arc.
  This proves the quasiconvexity. Now apply
the following result.

\begin{lemma} \label{Q-convex}
Suppose $E \subset {\Bbb R}^d$ has zero $d$-dimensional measure and
$\Omega={\Bbb R}^d \setminus E$  is quasiconvex.
Then any homeomorphism $f: {\Bbb R}^d \to {\Bbb R}^d$ which is locally
biLipschitz on $\Omega$ is biLipschitz on all of ${\Bbb R}^d$.
\end{lemma}

\begin{pf}
Given $x,y \in \Omega$, let $\gamma \subset \Omega$ be an arc
of length $\leq C|x-y|$ connecting them.
By integrating $|\nabla f|$  along $\gamma$ we see that $f$ is Lipschitz
on $\Omega$ and hence on all of ${\Bbb R}^d$. This means that
$f$ is absolutely continuous on lines and by the 
Radamacher-Stepanov theorem (e.g., Theorem 29.1 of \cite{Vaisala71}),
it is differentiable almost everywhere. Since $f$ is locally 
biLipschitz on a set of full measure, we deduce that 
$f$ is quasiconformal on ${\Bbb R}^d$ (the analytic definition 
of quasiconformality, Theorem 34.6 of \cite{Vaisala71}).
Thus  $f(E)$ has zero $d$-dimensional measure,
$f^{-1}$ is  globally quasiconformal and also
locally biLipschitz almost everywhere. Hence 
$F^{-1}$ is absolutely continuous on almost all lines.
This implies that given two points
$z,w$, we can connect them by a curve of length 
$\leq 2 |z-w|$ along which $f^{-1}$ is absolutely continuous
and has bounded derivative (just consider a family of connecting 
arcs which sweeps out positive measure). Integrating along the curve shows
$f^{-1}$ is also Lipschitz, as desired.
\end{pf}

To build a removable set 
which is weakly porus, we want to show that ordinary 
``flat'' cubes can be used in the previous construction, i.e, 
the ``flexible surfaces'' are not really needed for the quasiconformal
construction. We do this by proving that a flexible surface is 
actually a quasiconformal image of a flat square.

\begin{lemma}
Suppose $f$ is a diffeomorphism of
$V=[0,1]^2 \times [-\epsilon, \epsilon]$ into ${\Bbb R}^3$. 
Then there is a $\delta > 0$ and a homeomorphism 
$g: U \equiv [0,1]^2 \times [-\delta, \delta] \to {\Bbb R}^3$ such that
if $A = ([0, \delta] \cup [1-\delta,1])^2 \times \delta$, and 
$2A = ([0, 2\delta] \cup [1-2\delta,1])^2 \times \delta$, then 
\begin{enumerate}
\item $g=f$ on $A$.
\item $g$ is uniformly quasiconformal on   $U \setminus 2A$.
\item $g([0,1]^2)$ approximates $f([0,1]^2)$ to within $\delta$
in the Hausdorff metric.
\end{enumerate}
\end{lemma}

{\it Proof:}
Let $S$ be a flexible surface contained in $V$ and let 
$F$ be a quasiconformal approximation to $f$ on a
neighborhood $V$ of $S$.
We claim that there is a quasiconformal 
map $h$  of $U$ into $V \cup 2A$ which is the identity on $A$ and maps some 
subsquare of $[0,1]^2$ to the surface $S$. Given this,
define $g= F\circ h$. It easy to verify the desired properties, 
so we only have to construct the map $h$.

\begin{figure}[htbp]
\centerline{\psfig{figure=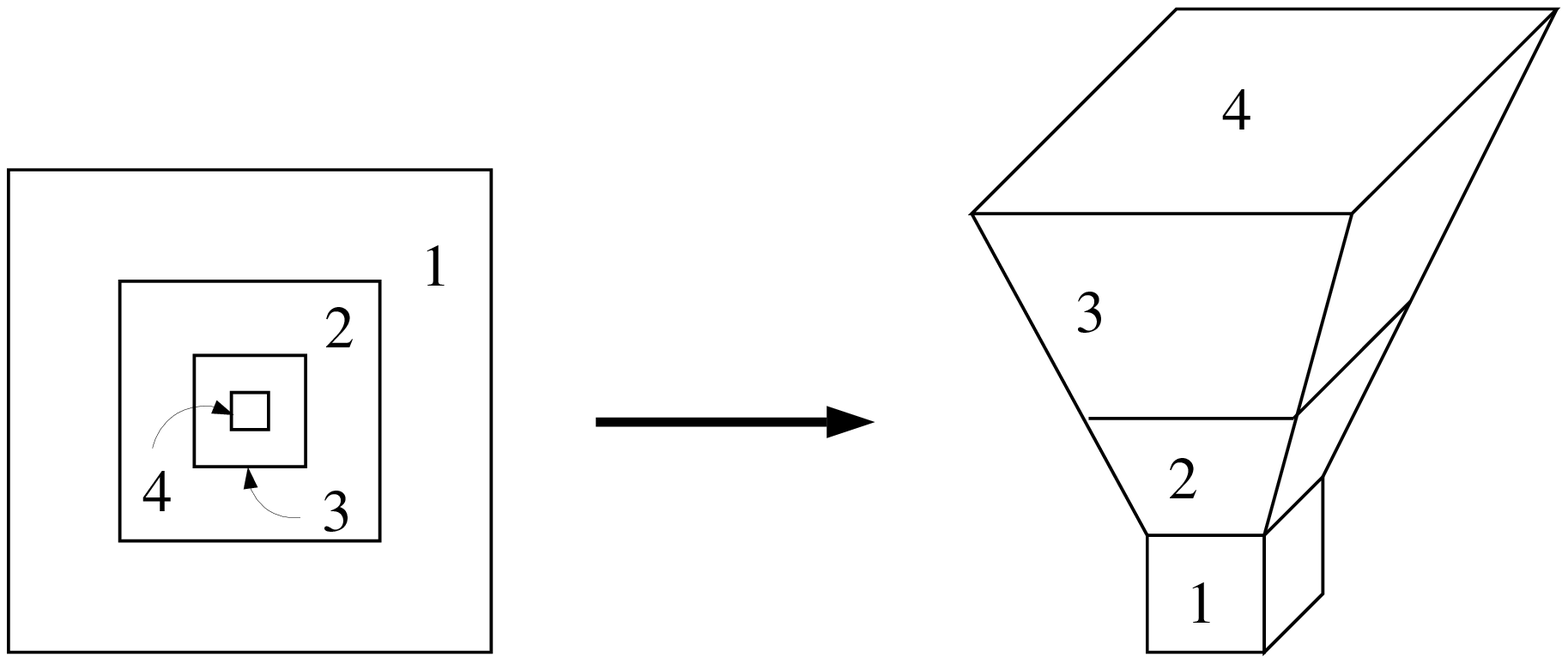,height=1.75in} }
\mycaption{\label{chimney} Step 1: 
       Quasiconformally map square to an expanding tower.} 
\end{figure}

\begin{figure}[htbp]
\centerline{ \psfig{figure=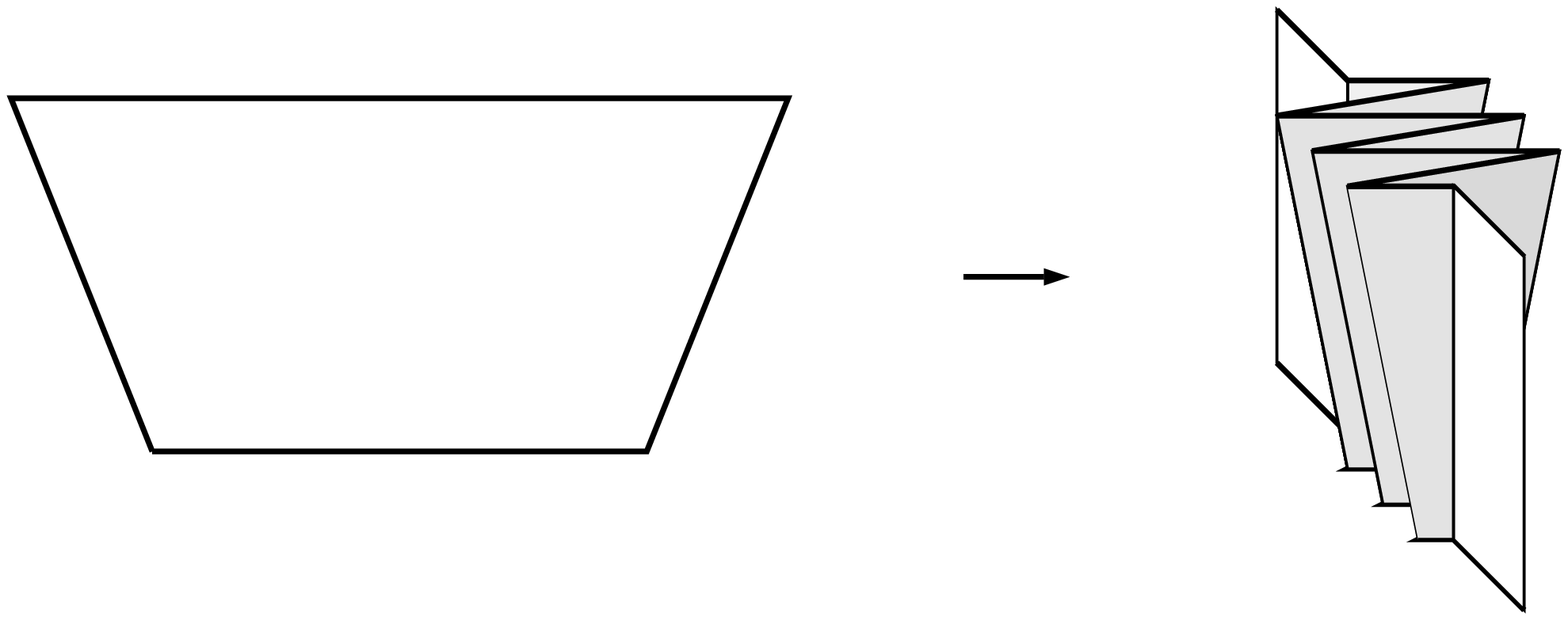,height=1.75in} }
 \mycaption{ \label{shrink-sides} Collapsing a side of expanding tower to a 
neighborhood of a side of a straight tower.} 
\end{figure}

\begin{figure}[htbp]
\centerline{ \psfig{figure=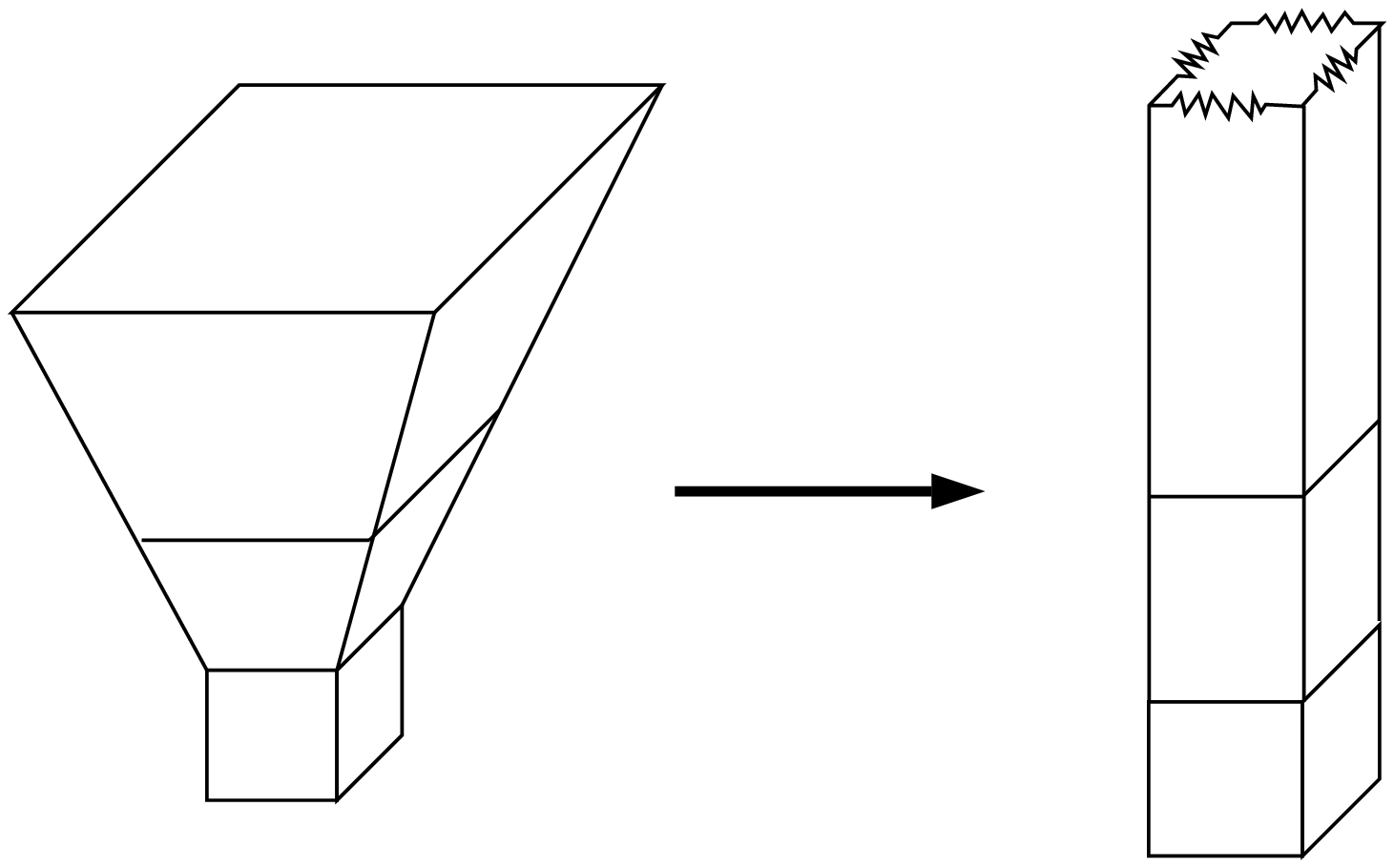,height=1.75in} }
 \mycaption{ \label{wiggly-chim} Step 2: Fold the expanding tower to a straight
 tower with a flexible surface on top. } 
\end{figure}

\begin{figure}[htbp]
\centerline{ \psfig{figure=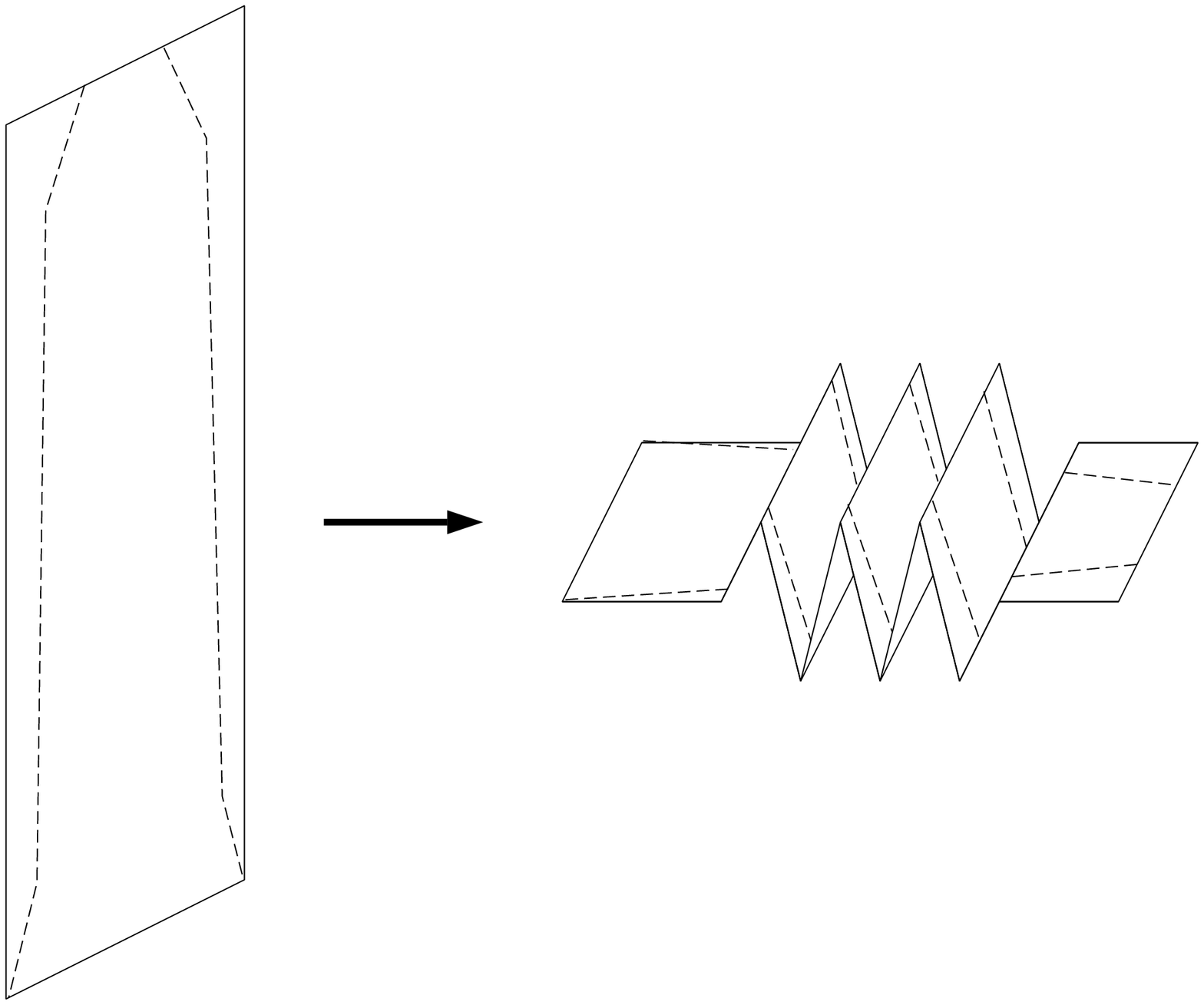,height=1.75in} }
 \mycaption{ \label{fold-str-side} Collapsing a straight side into 
the base of the tower.} 
\end{figure}

\begin{figure}[htbp]
\centerline{ \psfig{figure=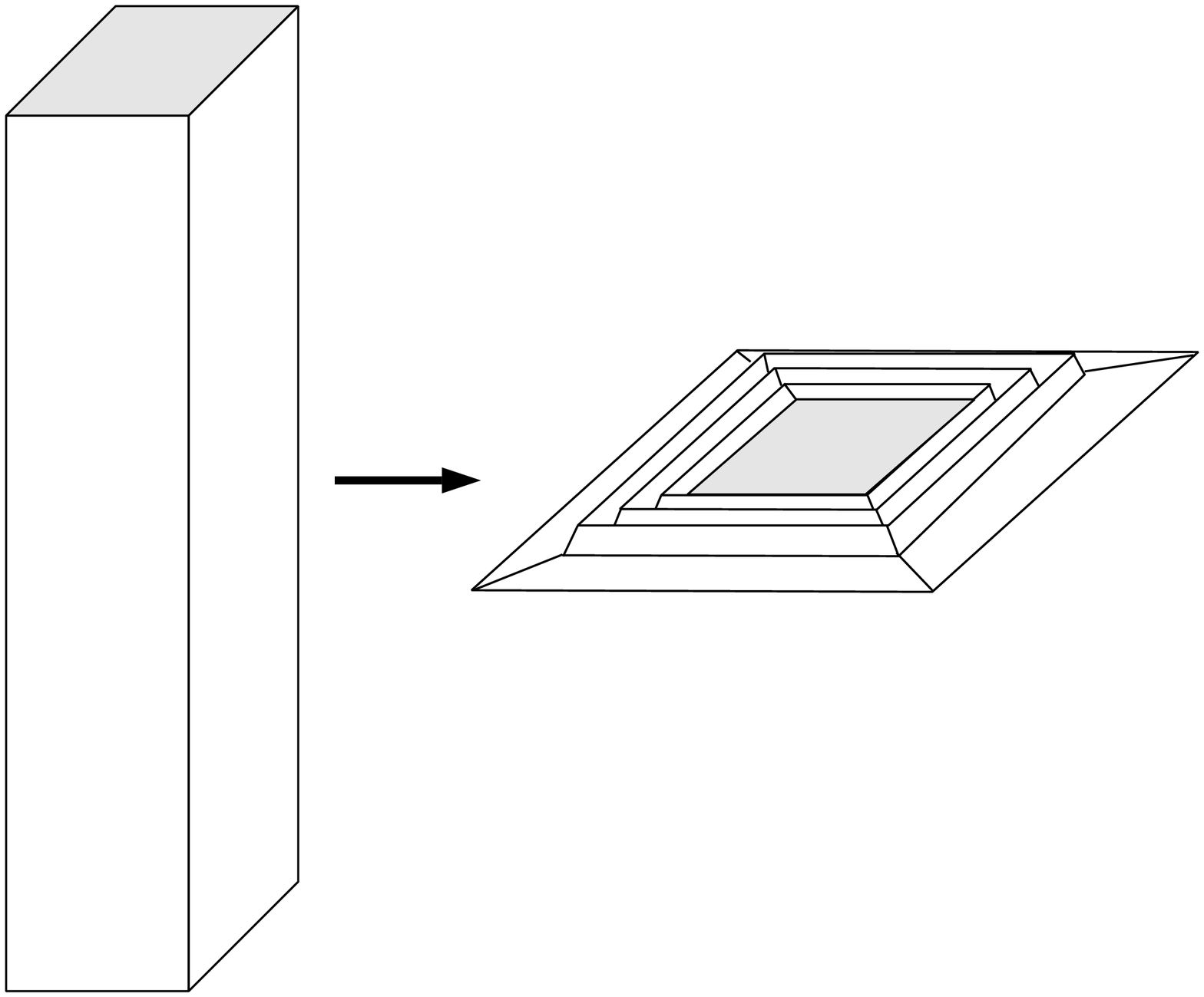,height=1.75in} }
 \mycaption{ \label{flatten-chim} Step 3:  Fold sides of straight tower 
to collapse it to a  neighborhood of base. } 
\end{figure}

This is easy to do in a couple  of steps. First, we can
quasiconformally  map the 
square to an ``expanding tower'' as in Figure \ref{chimney}.
The top of the tower is a large square which can be locally
biLipschitz mapped to a flexible square.
 The sides of the expanding tower can be 
folded as in Figure \ref{shrink-sides}.
to agree with the oscillation on the top. The result
is a surface which is close to a ``straight tower'', as in 
Figure \ref{wiggly-chim}.
 Finally, the sides of the straight tower can be folded 
as in Figure \ref{fold-str-side} to  ``collapse''
into a neighborhood of $[0,1]^2$, with the top mapping to 
$[\epsilon, 1-\epsilon]^2$. See Figure \ref{flatten-chim}.
(The straight side should  be mapped into the region 
bounded by the dotted line by a locally biLipschitz map
before folding; then after the folding the vertical projection
will be a trapezoid and the four sides will join together 
correctly).
  Composing these steps gives 
a uniformly quasiconformal map of a (very thin) neighborhood
of $[0,1]^2$ to a neighborhood of the  the desired surface,
This proves the lemma. $\qed$


The construction of non-removable sets now proceeds as before. The 
only difference is that instead of constructing quasiconformal 
approximations to arbitrary diffeomorphism, we now construct 
quasiconformal maps whose images approximate the images of the 
diffeomorphism (but the parameterizations do not necessary 
approximate each other).  However, this is sufficient.

Using the remarks above, we see that the set we construct can be made
disjoint from the faces of all dyadic cubes in ${\Bbb  }R^3$.  Thus the 
projection on each coordinate axis is totally disconnected which 
proves Corollary \ref{product-cor}.

\section{Small non-removable sets for biLipschitz maps in ${\Bbb R}^2$
and ${{\Bbb R}}^3$ }
\label{small-biLip-sec}

We now return to building non-removable sets for locally 
biLipschitz maps. In this section we show how to construct such 
sets with small Hausdorff measure. In the next section, we 
show how to insure that the image has small measure.

Since the proofs   in ${\Bbb R}^2$ and ${\Bbb R}^3$ are almost
identical, but easier to visualize in ${\Bbb R}^2$, we will 
consider that case first.
We will show that given a function 
$\varphi(t) = o(t)$, there is a totally disconnected 
$E \subset {\Bbb R}^2$ with $\Haus^\varphi(E)=0$
which is not removable for locally biLipschitz mappings.

Just as in Section \ref{small-QC-2}, we start with $E^0 = \{z: |z|=1\}$.
 Let $\Omega^0$ be the complement of $E^0$,   let $\Omega^0_0$ 
be its unbounded component and $ \Omega^0_1$ the 
bounded component.
Define 
$$ f^0_0(z) = z, \quad z \in \Omega^0_0,$$
$$ f^0_1(z) = \frac 12 z, \quad  z \in \Omega^0_1.$$

The induction hypothesis is as follows.
Suppose we are given a compact set $E^n$ which is a finite
union of $J = J_n$ smooth closed curves, $\{ E^n_j\}$. 
Let $\Omega^n$ be the complement
of $E^n$. Its unbounded component is denoted $\Omega^n_0$ and 
the bounded components are denoted $\Omega^n_j$, $j=1,\dots,J_n$.
Suppose we are given homeomorphisms $f^n_j$ on $\Omega^n_j$, $j=0,
\dots, J_n$, which are
locally biLipschitz  with constant $M$ on each component, and  
so that $f^n_j(\Omega^n_j)$ lies 
in  $Y^n_j$, the bounded complementary component of $f^n_0(E_j)$.

What follows is  a description of how to construct $f^{n+1}$
and $E^{n+1}$ from $f^n$ and $E^{n}$.

{\bf Step 1:} This is almost  exactly as in Section 
\ref{small-QC-2}, but with one small change.
As before, fix a  very small number $\eta_n>0$ 
and  for $j=1, \dots, J_n$ let
$U^n_j$ be an open topological annulus which  has   $E^n_j$ as
its ``outer'' boundary component and  which 
contains $  \{ z \in \Omega^n_j :  \dist(z, E^n_j) < \eta_n\}$.
Let $U^n =\cup_j U^n_j$.
 For $j=1, \dots, J_n$, let 
$$ \widetilde \Omega^n_j = \Omega^n_j \setminus \overline{U^n},
 \qquad  W^n_j  =  Y^n_j \setminus f^n_j (\widetilde \Omega^n_j)) .$$
Then $W^n = \cup_j W^n_j$ consists of $J_n$ annuli, so there is smooth 
diffeomorphism from  $U^n$ to $W^n$.
Thus we can construct a smooth diffeomorphism 
$g^n: {\Bbb R}^2 \to {\Bbb R}^2$ which agrees  with $f^n$ on 
${\Bbb R}^2 \setminus U^n$.

\begin{figure}[htbp]
\centerline{ \psfig{figure=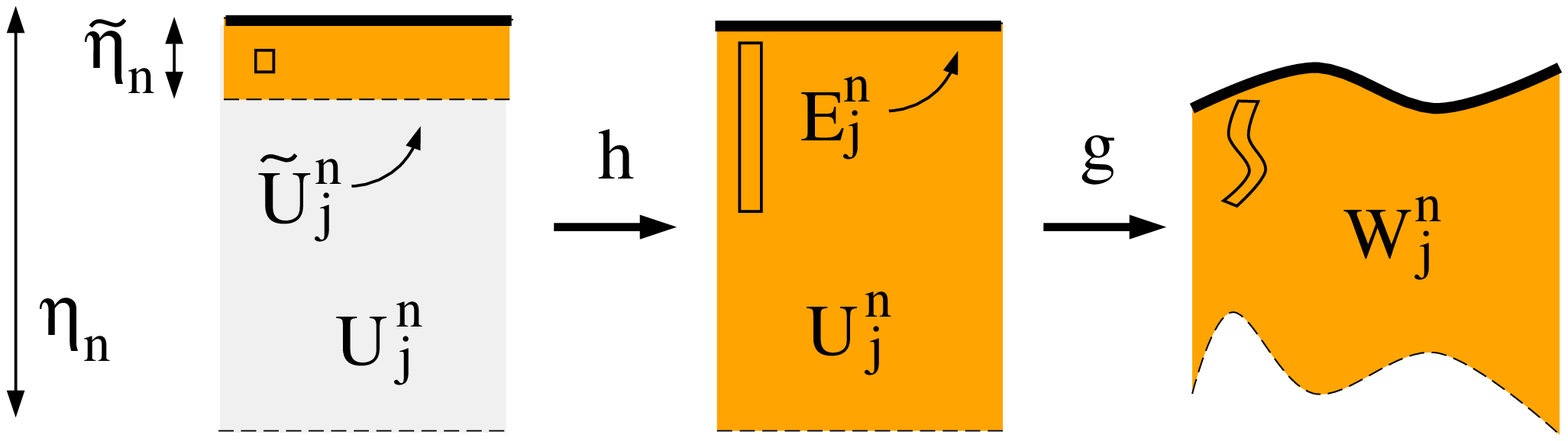,height=1.3in} }
 \mycaption{ \label{expand}  We may assume $g$ expands areas}
\end{figure}

In order for us to define $f^{n+1}_j$ to be biLipschitz later, it 
will be necessary to assume that $g^n$ is  area increasing on 
$U^n$.  To do this,  
replace $U^n$ by an even smaller neighborhood $ \widetilde U^n_j
\subset U^n_j$ of $E^n_j$, with 
$$ \widetilde U^n_j \subset \{ z: \dist (z, E^n_j) < \tilde  \eta_n << \eta_n
\} .$$
Now map $\widetilde U^n_j $ to $W^n_j$ by first taking the 
map $h: \widetilde U^n_j \to U^n_j$ which expands in the direction normal
to $E$ and then following with the map $g: U^n_j \to W^n_j$.
 The first map  expands  volume by a factor of $\eta_n/\tilde \eta_n
>>1$, so by selecting $\tilde \eta_n$ small enough (given the map $g$)
we can assume the composition is also area expanding.
See Figure \ref{expand}.
So replacing $\eta_n$ and $\tilde \eta_n$ and $U^n_j$ by $\tilde U^n_j$
if necessary,  we may assume we have a smooth mapping
$g^n: U^n_j \to W^n_j$ which expands area.

Now choose $\delta_n < \eta_n/10$ and consider the 
grid of $\delta_n \times \delta_n$ squares with 
vertices in $\delta_n {\Bbb Z} \times  \delta_n {\Bbb Z}$.
Let ${\cal S}_n$ be  a  collection of such squares which  cover 
$U^n$ and are contained in  $  \{ z :  \dist(z, E^n_j) < \eta_n/2\}$.
Let $F= \cup_{Q \in {\cal S}} \partial S$.  

{\bf Step 2:}  
As in Sections \ref{biLip-plane} and \ref{biLip-space}
 we replace the edges of the squares
by flexible arcs  to get a set $\widetilde F$, 
and we define $f^{n+1}$ on a neighborhood of 
this arcs to be a locally biLipschitz approximation to $g^n$ and
to agree with $f^n$ outside $U^n$.

We now have a smooth diffeomorphism $f^{n+1}$ defined on an open set
$\Omega^{n+1}_0$ which contains 
$\widetilde F \cup \cup_{j\geq 0}  \widetilde \Omega^n_{j}$. 
Without loss of generality we may assume that $\Omega^{n+1}_0$ is 
bounded by a finite number of smooth closed curves. Let 
$E^{n+1} = \partial \Omega^{n+1}_0$, and let $\Omega^{n+1}_j $, $j >1$ be an 
enumeration of the finitely many bounded complementary 
components. 

In Section \ref{small-QC-2},
 we  defined $f^{n+1}_j$ on $\Omega^{n+1}_j$ simply  
as a Euclidean  similarity which 
maps $\overline{\Omega^{n+1}_j}$ into $Y^{n+1}_j$.
However, this map might have to shrink the component a great 
deal to fit it inside $Y^{n+1}_j$, and we lose control of the 
biLipschitz constant.
However, because we have arranged for $g^n$ to be area 
expanding, we will be able to find a locally  biLipschitz mapping of 
a subdomain $\widehat \Omega \subset \Omega^{n+1}_j$ into $Y^{n+1}_j$.
First observe that
since $\Omega^{n+1}_j$ approximates a square $Q$ in $F$ as closely 
as we like, its area is as close as we like to the area 
of the squares in $F$. Similarly, $Y^{n+1}_j$ is an approximation 
to $g(Q)$, so we may assume its area is bigger than 
  $10 \text{area}(\Omega^{n+1}_j)$. 

\begin{figure}[htbp]
\centerline{ \psfig{figure=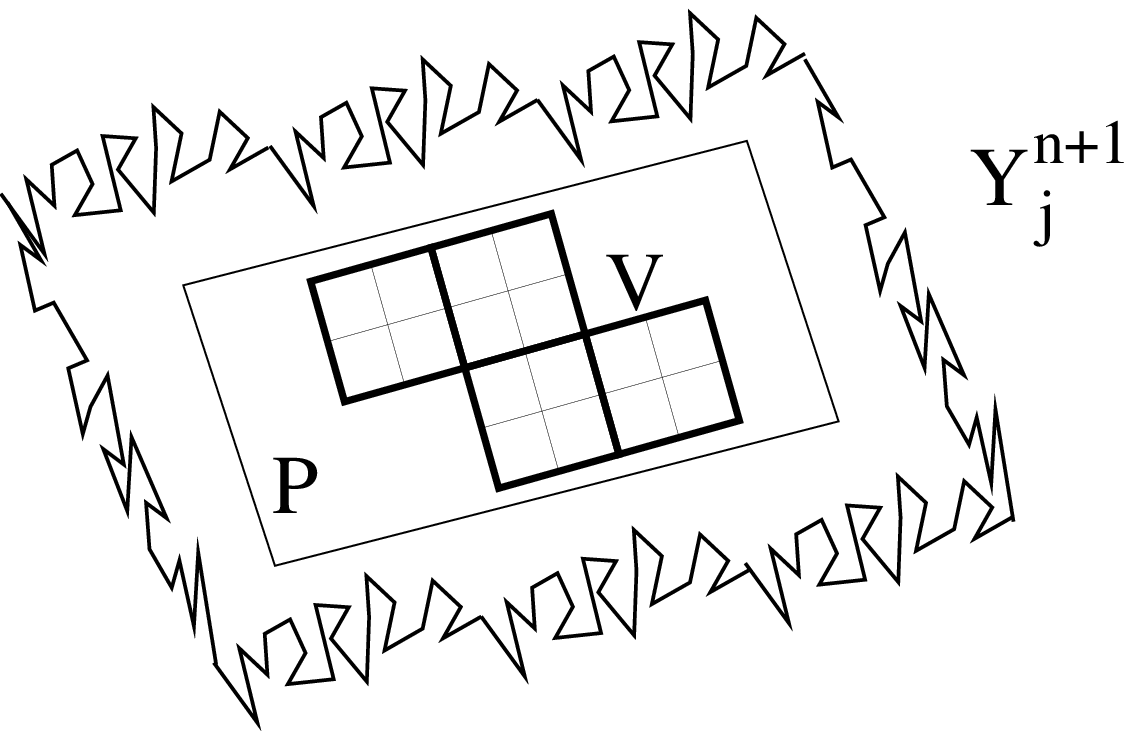,height=2.0in} }
 \mycaption{ \label{sqs-Y}   Squares in $Y^{n+1}_j$ } 
\end{figure}

Furthermore, since $g$ is smooth, if we 
take the squares in the construction small enough then $Y^{n+1}_j$
will approximate a parallelogram. Choose a true parallelogram 
$P \subset Y^{n+1}_j$ and choose a collection of disjoint squares
of size $\rho \times \rho$ 
in $P$ with connected union 
and which  cover at least half the area of $Y^{n+1}_j$ (and 
hence more than the area of $\Omega^{n+1}_j$). The number 
$\rho$ should be chosen so $\rho << \text{diam}(\Omega^{n+1}_j)$. See
Figure \ref{sqs-Y}.

\begin{figure}[htbp]
\centerline{ \psfig{figure=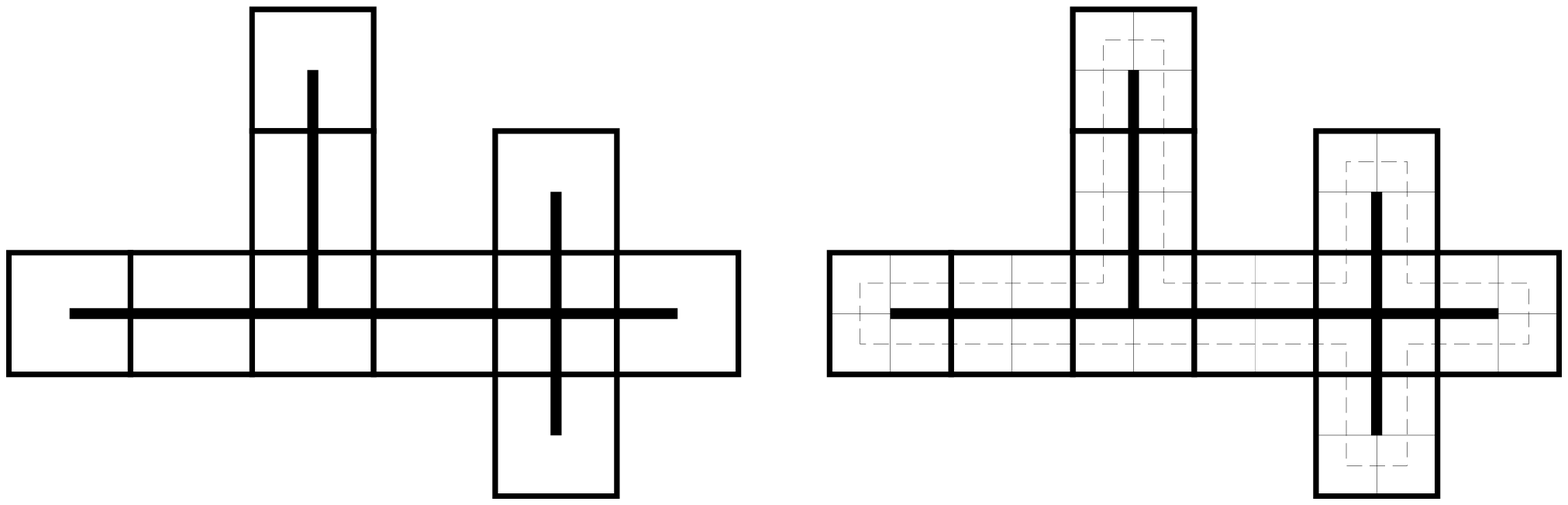,height=1.5in} }
 \mycaption{ \label{path}   Finding a Hamiltonian cycle in the 
``doubled''  graph of squares } 
\end{figure}

Since $Y^{n+1}_j$ is connected it is possible to choose the 
squares so that their union is connected. Thus we can think 
of the collection of squares as a graph where the squares are 
vertices and squares that share an edge are considered adjacent
in the graph. We want to label the squares $S_1, S_2, \dots$
so that $S_k$ and $S_{k+1}$ are adjacent, i.e., we want to 
find a Hamiltonian graph. Since the graph is connected, we 
can certainly find a spanning tree, but it may be impossible to find 
a Hamiltonian cycle. See the top picture in Figure \ref{path}.
However, if we replace each of our original squares by four squares
of half the size then it is always possible to find a 
Hamiltonian path. More  generally,
\begin{lemma}
Let ${\cal S} = \{Q_j\}$ be a connected collection of unit 
squares from the usual lattice in ${\Bbb R}^d$. Let 
${\cal S}' = \{ Q_{k}'\}$ be the collection obtained by replacing 
each cube in ${\cal S}$ by the $2^d$ subcubes of side length $\frac 12$.
Then the graph  $G$ with vertices ${\cal S}'$ and edges defined by 
adjacency of cube faces has a Hamiltonian cycle.
\end{lemma}

This is very easy by induction on the number of cubes, and we leave the 
proof to the reader.

\begin{figure}[htbp]
\centerline{ \psfig{figure=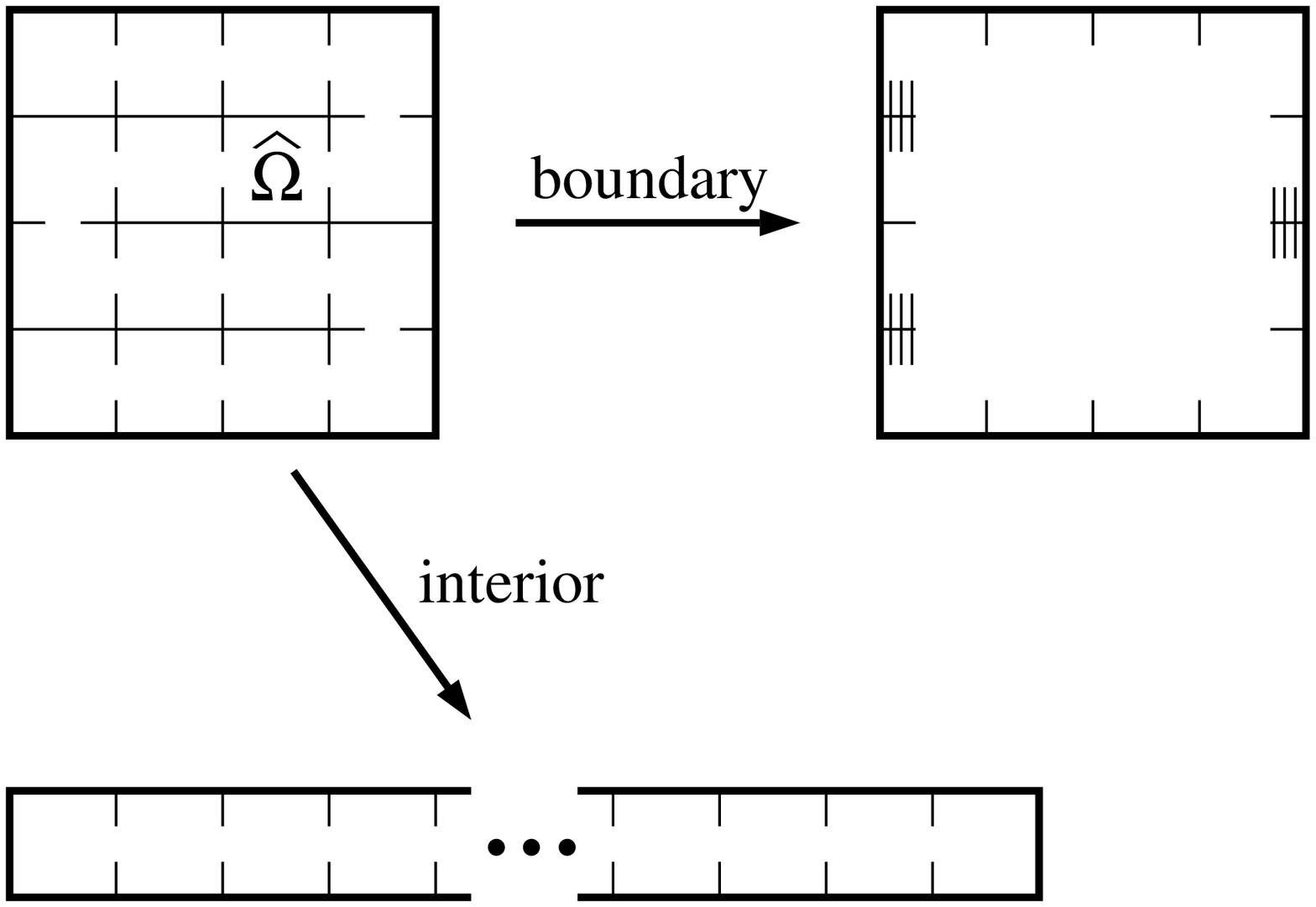,height=2.5in} }
 \mycaption{ \label{BL2block}  A flexible subdomain 
$\widehat \Omega$  of $\Omega^{n+1}_j$ and biLipschitz images
of its boundary and interior.}
\end{figure}

Define a subdomain $\widehat \Omega \subset \Omega^{n+1}_j$
as illustrated in the upper left of Figure \ref{BL2block}. This subdomain is 
topologically a disk, but looks like a decomposition of $\Omega^{n+1}_j$
into a chain of $\rho \times \rho$ squares.
Because $\Omega^{n+1}_j$ has much smaller area than $Y^{n+1}_j$,
the number of these squares is less than the number of squares chosen 
in $Y^{n+1}_j$ above.
If the arcs in $\partial \widehat \Omega$
are made up of flexible arcs which can be shrunk by locally biLipschitz
maps, then there is a locally biLipschitz map $h_1$ on a neighborhood
of $\partial \widehat \Omega$ which is the identity on $\partial 
\Omega^{n+1}_j$ and which maps $ \partial \widehat \Omega$
into $ g^{-1}(Y^{n+1}_j \setminus P)$.
 Thus 
$h \circ g$ is a locally biLipschitz mapping of $\partial \widehat 
\Omega$ into  $\partial Y^{n+1}_j \setminus P$.
By adjoining a neighborhood of $\partial \widehat 
\Omega$ to $\Omega^{n+1}_0$  for each $j$ we obtain a 
new region $\widetilde \Omega^{n+1}_0$ and an 
extension of $f^{n+1}_0$ to the new region.  The boundary of 
this expanded region is denoted $ \widetilde E^{n+1} = \cup 
\widetilde  E^{n+1}_j$, each 
component of which we may assume to be a smooth closed curve.

\begin{figure}[htbp]
\centerline{ \psfig{figure=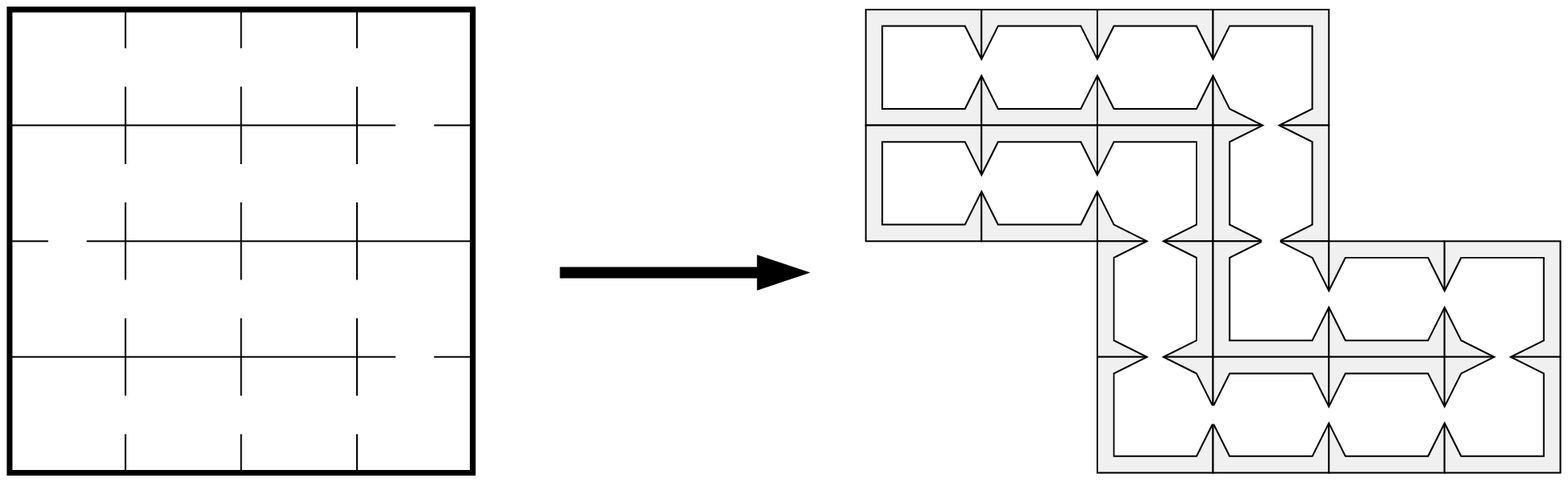,height=1.5in} }
 \mycaption{ \label{BL2int}  The  biLipschitz image of $\widehat \Omega$
inside $Y^{n+1}_j$. }
\end{figure}

On the other hand, is easy to see
that $\widehat \Omega$ itself  can be  locally biLipschitz  mapped to  
the long narrow region in the bottom of Figure \ref{BL2block}.
This in turn can be  locally biLipschitz mapped into the any region which 
is a ``chain '' of  similar number of similarly sized squares.
In particular, it can be mapped to the squares in $P$. With a
slight adjustment we can easily make the image a Jordan domain, 
e.g., see Figure \ref{BL2int}.
Thus on each bounded complementary component of $ \widetilde E^{n+1}$ we have a 
uniformly biLipschitz mapping into the bounded complementary 
component of $f^{n+1}_0(\widetilde E^{n+1}_0)$.
 This completes the proof of the inductive step.

Passing to a limit exactly as before we obtain a homeomorphism 
$f$ of ${\Bbb R}^2$ which is   uniformly locally biLipschitz  of
a totally disconnected set $E$. Since the construction shows
that we can take 
$$E \subset \{ z: \dist(z, E_n) \leq \eta_n\},$$
where $E_n$ is finite union of smooth curves and $\eta_n$ 
is as small as we wish (independent of $E_n$) it is easy to 
construct $E$ so $\Haus^\varphi(E) =0$.

We now  make a few comments on how to 
modify the construction of the  so that it
works in ${\Bbb R}^3$. Just as above we may assume 
we have a mapping $g: \Omega^{n+1}_j \to Y^{n+1}_j$ which is 
volume expanding and almost linear. Moreover, $\Omega^{n+1}_j$
is a close approximation to a cube and $Y^{n+1}_j$ is a 
close approximation to a parallelepiped. We want to construct
a subdomain $\widehat \Omega \subset \Omega^{n+1}_j$ so that
$\Omega^{n+1}_j \setminus \widehat \Omega$ is a union of 
surfaces which can be mapped to a given neighborhood of 
$\partial \Omega^{n+1}_j$ by locally biLipschitz mapping,
and so that there is a locally biLipschitz
mapping from $\widehat \Omega $ to $Y^{n+1}_j$ which misses  a
given neighborhood of $\partial Y^{n+1}_j$.
As before we find a real parallelepiped $ P \subset Y^{n+1}_j$
and a collection of cubes in $P$ which cover half the volume
of $Y^{n+1}_j$. As in the previous section, it may not be true 
that there is a Hamiltonian path in the resulting graph of 
cubes, but if we replace each cube by the 8 subcubes of half
the size the resulting graph 
always has a Hamiltonian cycle.

As before, we construct $\widehat \Omega$ by 
 dividing $\Omega^{n+1}_j$ into small ``cubes''
of approximately the same size as those chosen in $P$ above.
Now replace the flat sides of the cubes by copies of
``flexible surfaces'' constructed in Section \ref{flex-square-sec}.
The flexible surfaces are chosen so that there is biLipschitz 
mapping of a neighborhood of the union of faces  into 
$g^{-1}(Y^{n+1} \setminus P)$,  which is a neighborhood of 
$\partial \Omega^{n+1}_j$. Add this neighborhood to
$\Omega^{n+1}_0$ and let $\{\widehat  \Omega_k\} $ be 
the complement in $\Omega^n_j$.

So far $\widehat \Omega$ is a disjoint union of many topological balls.
To make it a single ball, enumerate the components  $\widehat
\Omega_k$ so that $\widehat \Omega_k$ and $\widehat \Omega_{k+1}$
share a face and join adjacent cubes.
This makes gives us a single connected component  $\Omega^{n+1}_j$
which is topologically a ball and which can be biLipschitz mapped 
to a chain of cubes as in the previous section.
This in turn can 
be biLipschitz mapped into $P \subset Y^{n+1}_j$, just as in the 
previous section.

We now have a biLipschitz map $f^{n+1}_0$ defined on the 
open set $\Omega^{n+1}_0$ with boundary $E^{n+1} = \cup_j 
E^{n+1}_j$, so that each $E^{n+1}_j$ is diffeomorphic to the
$2$-sphere and bounds a topological $3$-ball $\Omega^{n+1}_j$.
Moreover there is a uniformly biLipschitz map of $\Omega^{n+1}_j$
into $Y^{n+1}_j$, the bounded  complementary component of 
$f^{n+1}_0 (E^{n+1}_j)$.  This completes the induction step.

The proof that in the limit we get  homeomorphism $f$ which 
is uniformly locally  biLipschitz  off a Cantor set $E$, is 
just as before. Similarly for the proof  that given 
 $\varphi(t)= o(t^2)$ we may construct $E$ so that 
$\Haus^\varphi(E)=0$.

\section{Making $f(E)$ small}
\label{biLip-small-image}

In this section we will show how to modify the construction 
in order to insure $\Haus^\varphi(f(E))=0$. We start
by reviewing two additional facts we will use.

The first is  a result of 
Dacorogna and Moser \cite{Dacorogna-Moser} that if 
$f: \Omega_1 \to \Omega_2$ is a diffeomorphism of smooth 
domains of equal volume, then there is another diffeomorphism 
$f_0: \Omega_1 \to \Omega_2$ which agrees with $f$ on 
$\partial \Omega_1$, but which is volume preserving.
(In fact, we can specify the Jacobian any smooth way 
we want as long the total volumes work out correctly).
The second fact is a generalization of the argument 
used in Section \ref{small-QC-space}.

\begin{lemma} \label{bi-sub}
 Suppose $\Omega_1$ and $\Omega_2$ are two smooth domains in 
${\Bbb R}^d$ with finite volume. Then there are subdomains 
$\widehat \Omega_i \subset \Omega_i$, $i=1,2$ so that 
$\Omega_i \setminus \widehat \Omega_i$, $i=1,2$  has finite
$(d-1)$ dimensional measure and such that there is a 
locally biLipschitz map $f: \widehat \Omega_1 \to \widehat
\Omega_2$ (and the biLipschitz constant depends only on
the ratio of the volumes of $\Omega_1, \Omega_2$).
\end{lemma}

\begin{pf}
First note that each smooth domain  can be 
biLipschitz mapped to a domain which is union of cubes
of side length $\rho$ which has comparable volume to the 
original domain. See the top of Figure \ref{bl-equal}.
\begin{figure}[htbp]
\centerline{ \psfig{figure=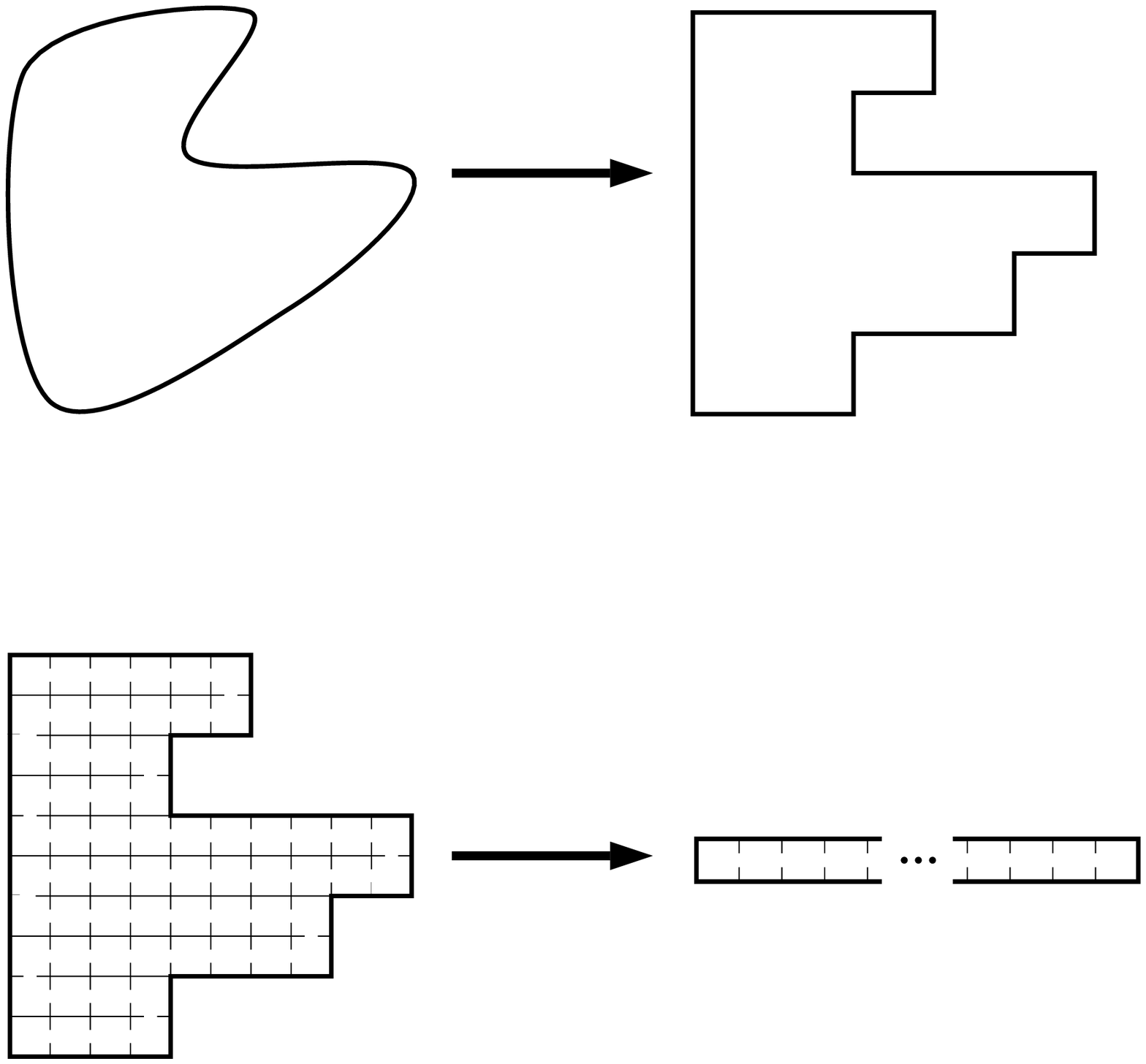,height=2.5in} }
 \mycaption{ \label{bl-equal}  A smooth domain can be biLipschitz
mapped to a union of cubes } 
\end{figure}
By replacing the cubes by cubes of half the size one can 
insure a Hamiltonian path in the resulting graph. 
Now define the  subdomain by creating small openings
(say of size $\rho/10$)
between adjacent squares  along the Hamiltonian paths.
See the bottom of Figure \ref{bl-equal}.
The resulting domains are clearly locally biLipschitz 
equivalent to tube of width $\rho$ and the correct volume, 
and hence are locally  biLipschitz equivalent with each other.
\end{pf}

\begin{figure}[htbp]
\centerline{ \psfig{figure=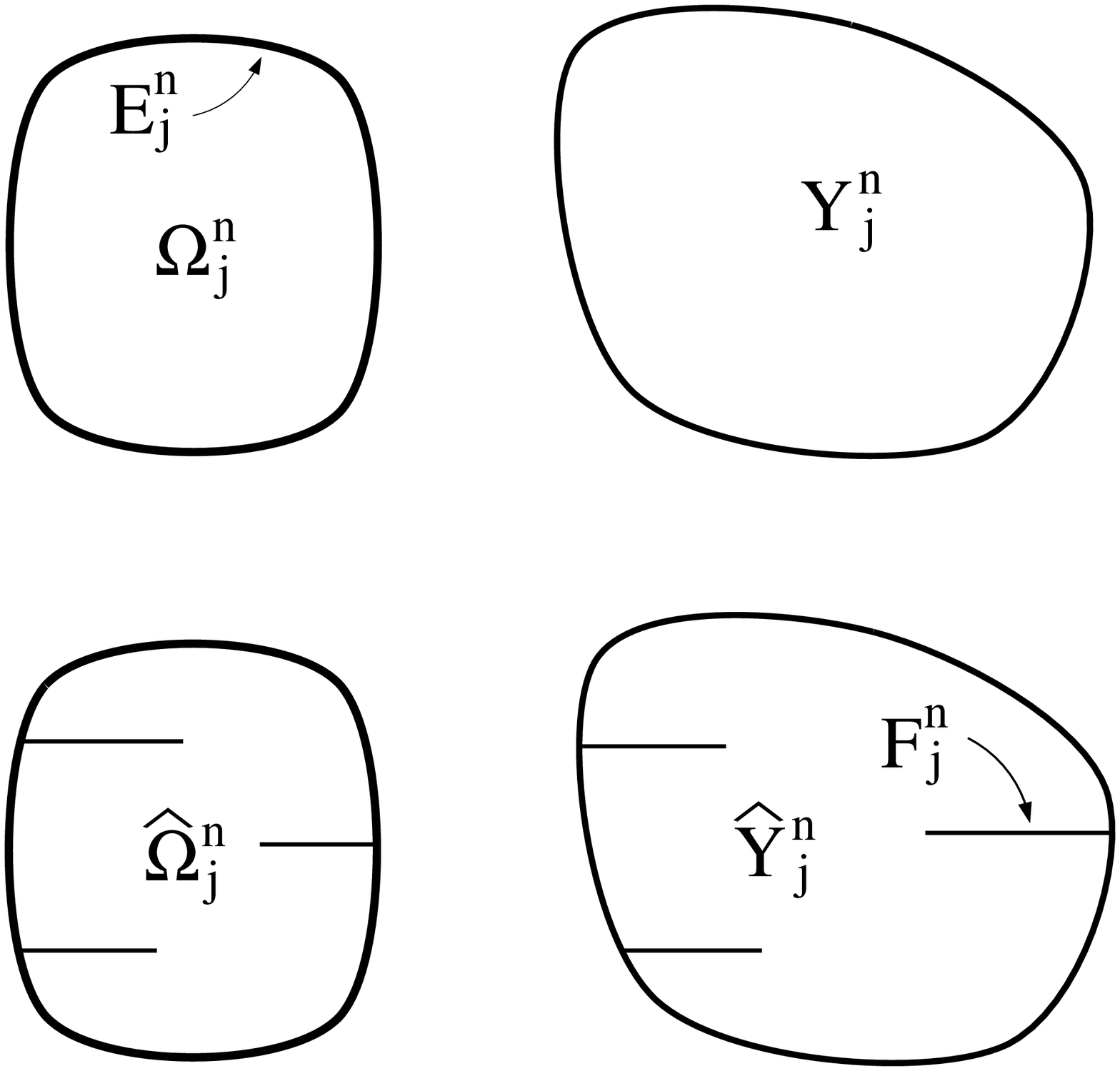,height=2.5in} }
 \mycaption{ \label{defns1}  Recalling the definitions } 
\end{figure}

We can now  make the desired modifications of the construction.
The new part of the induction hypothesis is 
that we have a locally  biLipschitz map $f^n_0$ on $\Omega^n_0$ 
and locally biLipschitz maps   $f^n_j$ on subdomains  $\widehat 
\Omega^n_j \subset \Omega^n_j$, $j=1, \dots, J_n$.  
If  $F^n_j = \overline{ Y^n_j }\setminus  \widehat Y^n_j$
(where $ \widehat Y^n_j = f(\widehat \Omega_j^n)$),
then we also assume that  $F^n_j$ 
has finite $d-1$ dimensional measure and $\Omega^n_j 
\setminus \widehat \Omega^n_j$ is a finite union of 
flat surfaces (in fact is a union of faces of dyadic cubes).

 We cover 
both $ \partial \widehat \Omega^n_j$ and $F^n = \cup_j F^n_j$
 by balls so that the $\varphi$-sum 
of the radii is small (say less than $\frac 1n$).
%
%
Next choose a topological annulus $W^n_j \subset Y^n_j$
which covers $F^n_j$, whose outer boundary
is $\partial Y^n_j$ and which 
is contained in the  good covering of $F^n_j$ described above.
Then $ D^n_j =  \widehat \Omega^n_j \setminus (f^n_j)^{-1} (W^n_j)$
 is a Jordan subdomain of $\widehat \Omega^n_j$ and
$\widehat \Omega^n_j \setminus D^n_j$ 
has volume comparable to $W^n_j$. By changing $f^n_j$
slightly we can assume the volumes are equal.
\begin{figure}[htbp]
\centerline{ \psfig{figure=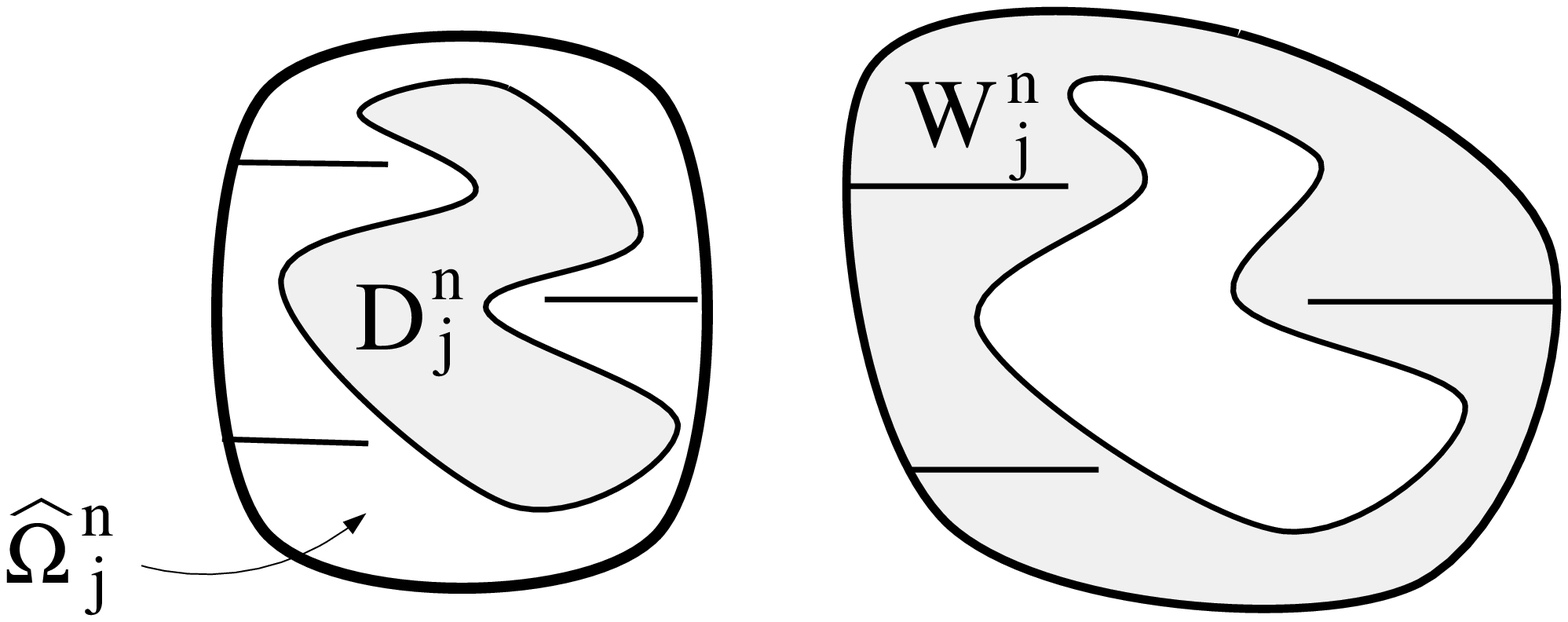,height=1.25in} }
 \mycaption{ \label{defns2}  Definition of $D^n_j$ and 
$W^n_j$ } 
\end{figure}

Recall that $\partial \widehat \Omega^n_j$ is a finite 
union of flat surfaces, so we divide it up into 
small squares and replace each by a flexible square.
The resulting domain is called $\widetilde \Omega^n_j$,
We choose the flexible surfaces   
so the map $f^n_0$ can be extended
to a uniformly biLipschitz map on $\Omega^n_0 \cup 
E^n_j \cup \partial \widetilde \Omega^n_j$ which maps 
a neighborhood of 
$  \partial \widetilde \Omega^n_j \cap \Omega^n_j$ into $W^n_j $.

\begin{figure}[htbp]
\centerline{ \psfig{figure=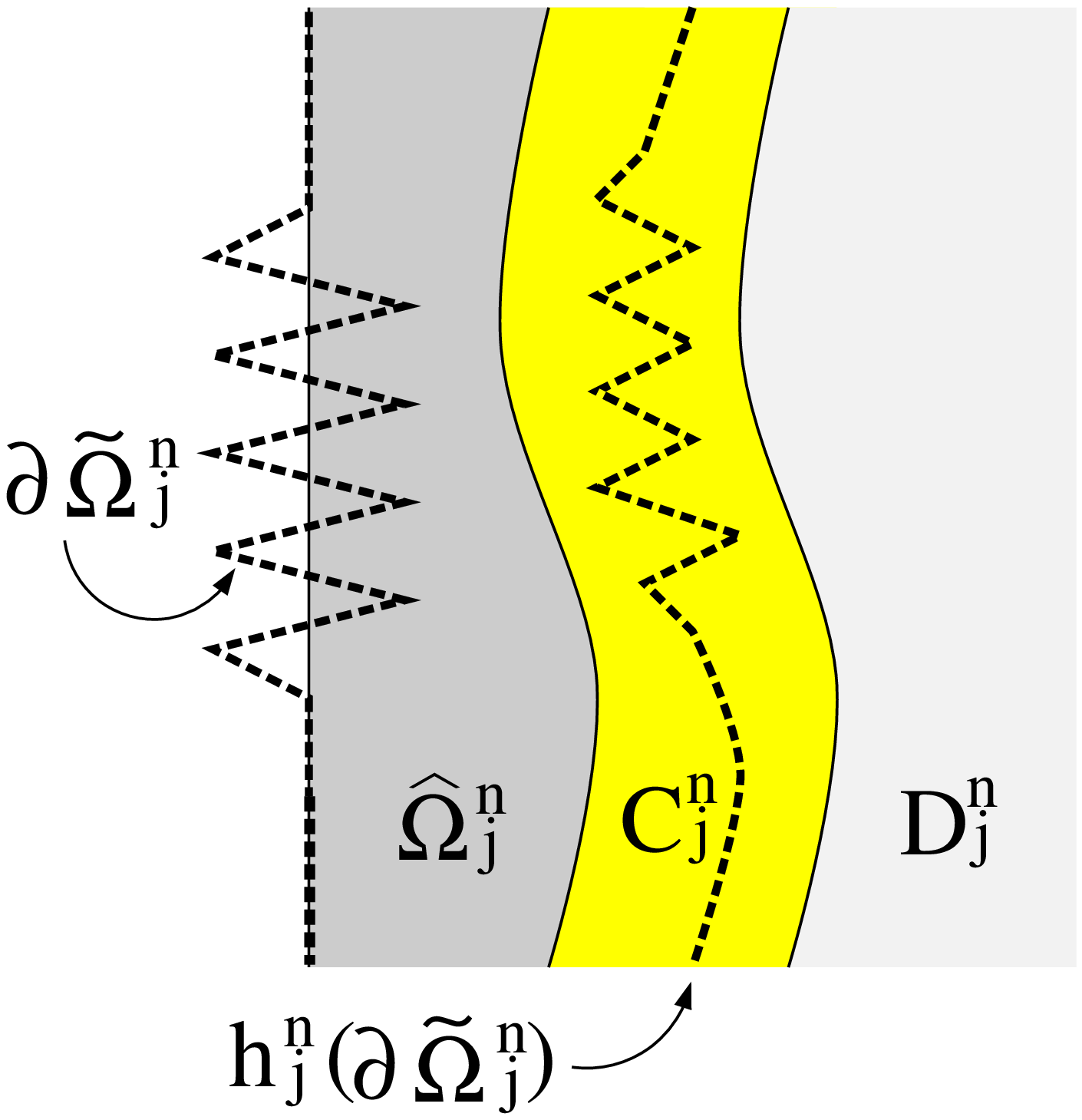,height=2.0in} }
 \mycaption{ \label{defnsC}  Definition of $h^n_j$  }
\end{figure}

Also, by choosing the flexible squares to be close
enough to the surfaces they replace, we may assume that
$\widetilde \Omega^n_j$ compactly contains $D^n_j$.
Moreover,  we may assume there is uniformly locally biLipschitz mapping
$h^n_j$  from 
$\widetilde \Omega^n_j$ into a neighborhood  $C^n_j$ of 
$\overline{D^n_j}$ which is the identity on  $D^n_j$.
We may also assume  $\overline{C^n_j} \subset \widehat \Omega^n_j$.
See Figure \ref{defnsC}.

\begin{figure}[htbp]
\centerline{ \psfig{figure=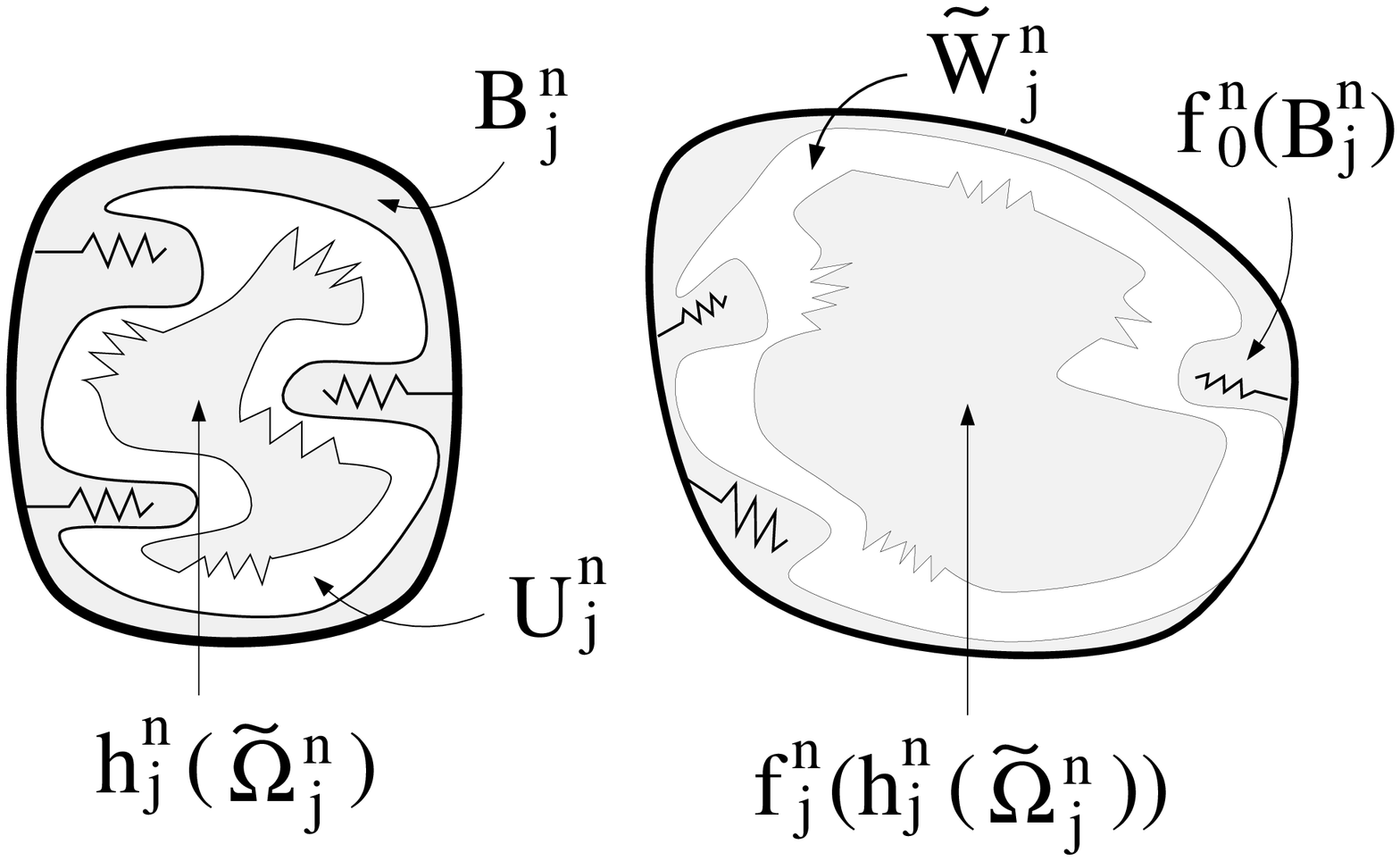,height=2.0in} }
 \mycaption{ \label{defns3}  Definition of $B^n_j$ and  $\widetilde W^n_j$ } 
\end{figure}

Now  define a topological annulus $B^n_j$ so 
that (see Figure \ref{defns3})  \newline
(1) $E^n_j$ is its outer boundary,  \newline
(2) it covers $\partial \widetilde \Omega^n_j \cap \Omega^n_j$, \newline
(3) it is disjoint from $C^n_j$ and  \newline
(4) the mapping $f^n_0$ which we extended to $\partial 
\widetilde \Omega^n_j$ has a uniformly locally  biLipschitz extension 
to $B^n_j$. 

We may take the volume of  $B^n_j $
to be as small as we like, in particular, smaller 
than   half the volume of $\widehat \Omega^n_j \setminus 
C^n_j$. Thus the volume of the annulus
 $$ U^n_j =\widehat \Omega^n_j \setminus (B^n_j \cup h^n_j(\widetilde
\Omega^n_j))$$
 is 
comparable to the volume of the annulus 
$$ \widetilde W^n_j = 
W^n_j \setminus (f^n_0(B^n_j) \cup f^n_j \circ h^n_j(\widetilde
\Omega^n_j)).$$

Then $U^n_j$ and $W^n_j$ are annuli of comparable 
volume (independent of $n$) and we have a 
locally biLipschitz map  $f^n_0$ from 
the outer boundary  component (i.e., $E^n_j$) of
$U^n_j$ to the outer boundary 
of $\widetilde W^n_j$ and we have  locally biLipschitz 
maps  $f^n_j \circ h^n_j$  of the bounded complementary component of each 
$U^n_j$ to the corresponding bounded complementary component
of $\widetilde W^n_j$.

 We now  use
the  result of Dacorogna and Moser \cite {Dacorogna-Moser} 
described at the beginning of this section
 to find diffeomorphisms $g_j: U^n_j \to 
W^n_j$ which extend $f^n_0$ and $f^n_j$  on the two boundary 
components and which 
multiply volumes by a constant factor (the ratio
of the volumes).

Now  proceed as before,  covering $U^n$ by cubes, 
replacing the faces by flexible surfaces, approximating 
$g_j$ on these surfaces and getting in the end  components $\Omega^{n+1}_j$
and $Y^{n+1}_j$ which have comparable volumes.
For each component $\Omega^{n+1}_j$ we  use Lemma 
\ref{bi-sub} to define a subdomain 
$\widehat \Omega_j$ which can be mapped into 
a subdomain of $Y^{n+1}_j$.

This completes the induction. Since we began the inductive step 
by insuring that our construction took place within 
a good covering, 
it is easy to see that the limiting set $E$ and homeomorphism 
$f$ satisfy $\Haus^\varphi(E) = \Haus^\varphi(f(E)) =0$.

As a final remark we observe that Corollary \ref{app-cor} is 
almost immediate.
If $f: \Omega_1 \to \Omega_2$ is a diffeomorphism, then we can write
$\Omega_1$ as a union of cubes so that $f$ is close to linear on
each cube. We replace the faces of these cubes by flexible surfaces 
and apply the construction and we obtain the desired homeomorphism.
Because the diffeomorphism may change volumes, we can only get a 
quasiconformal approximation.
If $\Omega_1$ and $\Omega_2$ are diffeomorphic by a volume preserving 
map then the construction of the last two sections applies and 
we can get a locally bi-Lipschitz approximation.


\bibliography{refs}
\bibliographystyle{plain}

\end{document}